\documentclass{article}

\usepackage{arxiv}
\usepackage{textcomp}

\usepackage{graphicx}
\usepackage{amsmath}
\usepackage[version=4]{mhchem}
\usepackage{siunitx}
\usepackage{longtable,tabularx}
\usepackage{url}
\usepackage{amsmath}
\usepackage{bm}
\usepackage{enumitem}
\usepackage{float}
\usepackage{natbib}
\usepackage{caption}
\usepackage{subcaption}
\usepackage{array}
\usepackage{amsfonts}
\usepackage[dvipsnames]{xcolor}
\usepackage{booktabs}
\usepackage{nomencl}
\makenomenclature
\usepackage{soul}
\usepackage{arydshln}
\usepackage{placeins}
\usepackage{multirow}

\usepackage{algorithm}
\usepackage[noend]{algpseudocode}

\usepackage{geometry}
\setcitestyle{numbers} 
\begin{document}
\title{SeAr PC: Sensitivity Enhanced Arbitrary Polynomial Chaos}


\author{
    Nick Pepper \\
	The Alan Turing Institute\\
	The British Library\\
	London, UK \\
	\texttt{npepper@turing.ac.uk} \\
    \And
    Francesco Montomoli  \\
    Department of Aeronautics \\
	Imperial College London\\
    London, UK
     \And
    Kyriakos Kantarakias \\
    Department of Aeronautics \\
	Imperial College London\\
    London, UK 
}

\maketitle

\begin{abstract}
This paper presents a method for performing Uncertainty Quantification in high-dimensional uncertain spaces {by combining arbitrary polynomial chaos with a recently proposed scheme for sensitivity enhancement} \cite{se-pce}. Including available sensitivity information offers a way to mitigate the \emph{curse of dimensionality} in Polynomial Chaos Expansions (PCEs). Coupling the sensitivity enhancement to arbitrary Polynomial Chaos allows the formulation to be extended to a wide range of stochastic processes, including multi-modal, fat-tailed, and truncated probability distributions. In so doing, this work addresses two of the barriers to widespread industrial application of PCEs. The method is demonstrated for a number of synthetic test cases, including an uncertainty analysis of a Finite Element structure, determined using Topology Optimisation, with 306 uncertain inputs. We demonstrate that by exploiting sensitivity information, PCEs can feasibly be applied to such problems and through the Sobol sensitivity indices, can allow a designer to easily visualise the spatial distribution of the contributions to uncertainty in the structure. 
\end{abstract}

\section{Introduction}
In any engineering product, some properties or environmental conditions will be uncertain. The objective of robust design is to achieve a final design that offers good performance across a range of uncertain conditions. Uncertainty Quantification (UQ) is used to evaluate the effect of uncertain parameters on a system, where we denote $\boldsymbol{\xi}\in \Xi \subseteq \Re^{n_u}$ the $n_u$ uncertain parameters; with $\mathcal{M}(\boldsymbol{\xi})$ the system response. Given a joint density $f(\boldsymbol{\xi})$, the designer is typically interested in estimating the mean and standard deviation of the response \cite{ZANG2005315}: 

\begin{align}
    &\mathbb{E}(\mathcal{M})=\int \mathcal{M}(\boldsymbol{\xi})f(\boldsymbol{\xi})d\boldsymbol{\xi},     \label{mean}\\
    &\text{Var}(\mathcal{M})=\int (\mathcal{M}(\boldsymbol{\xi})-\mathbb{E}(\mathcal{M}))^2f(\boldsymbol{\xi})d\boldsymbol{\xi}. \nonumber
\end{align}
Should $\mathcal{M}$ be inexpensive to compute, then the statistical moments of $f(\mathcal{M})$ can be readily calculated through Monte Carlo sampling. However, in most cases in Engineering $\mathcal{M}$ is expensive, with an evaluation requiring a single run of an expensive Computational Fluid Dynamics (CFD) \cite{cfd} or Finite Element (FE) simulation \cite{BLATMAN2010183}. In such cases a surrogate for $\mathcal{M}$ is required. Of the range of surrogate models proposed in the literature, one of the most popular has been the Polynomial Chaos Expansion (PCE), in which the uncertain Quantitiy of Interest (QoI) is expressed as the weighted sum of a multivariate polynomial basis. First proposed by Wiener \cite{homoChaos} for Gaussian stochastic processes, this result was later extended by Xiu and Karniadakis \cite{xiu} for other types of stochastic process. For certain stochastic processes, polynomial bases within the Askey scheme demonstrated optimal convergence. 

PCEs have been applied in a number of fields, with examples ranging from CFD \cite{CHATZIMANOLAKIS2019207,Kantarakias1}, medicine \cite{med}, ecology \cite{Pepper}, chaotic dynamics \cite{PREUQKantar,a13040090,Kantarakias_Shawki_Papadakis_2020} electronics \cite{electric},  and even astrophysics \cite{astro}. For a PCE of order $p$ there are $P+1$ linear combination terms in the expansion, calculated as:
\begin{align}
    P+1=\frac{(n_u+p)!}{n_u!p!}.
    \label{cardinality}
\end{align}
This can be problematic as $P$ scales rapidly with $n_u$. $P+1$ represents a necessary condition on the number of evaluations of $\mathcal{M}$ required to determine a unique solution for the coefficients of the expansion through a least squares (LSQ) fitting. However, in practice an oversampling will be necessary to achieve accurate results, perhaps requiring $2P-3P$ evaluations to estimate the coefficients. In any case, the number of required model evaluations greatly increases as $n_u$ is increased, a phenomenon referred to as the \emph{curse of dimensionality} \cite{eldred2009comparison}. This is a major obstacle to the massive industrial application of PCE. As an example, if a designer wishes to consider the effect of geometric uncertainties, then $n_u$ can grow to tens, if not hundreds of uncertain parameters (see, e.g. \cite{TO, turbo}). Accounting for uncertainties in loading, boundary conditions, and material properties, in addition to the geometry, grows the number of uncertain parameters further still \cite{KESHAVARZZADEH2017120}. At present, performing a Topology Optimisation while considering parametric uncertainty is only possible by reducing the dimensions of the uncertain space, for instance through Karhunen-Loeve Expansions (KLE) \cite{chen2010level}, or through reduced order models of the system \cite{MAUTE2009450}. Such methods can make a robust design optimization computationally tractable, but limits the types of stochastic processes that can be modelled, for instance to Gaussian random fields \cite{lazarov2012topology}. Similarly, the description of the parametric uncertainty can be made less rich by using intervals, rather than the joint density, in order to make UQ computationally tractable (see, e.g. \cite{WU201636}).


Several solutions have been proposed in the Polynomial Chaos literature for mitigating the \emph{curse of dimensionality}. For instance Smolyak's algorithm can be applied to assemble a sparse sampling grid on which the evaluate $\mathcal{M}$ \cite{smolyak1, smolyak2, smolyak3}. Another family of methods employ sparse regression techniques such as Least Angle Regression from compressive sensing in order to define sparse PCEs. By controlling the interaction order, that defines the coupling between various uncertain inputs in the PCE, the number of terms required can be reduced, with an adaptive scheme used to find the optimal polynomial basis \cite{L_then_2021}. However, as will be seen, these methods can still require more evaluations of $\mathcal{M}$ than are feasible when scaled to uncertain spaces with tens or hundreds of dimensions.   

A further challenge to the industrial application of PCE is in the handling of datasets that may contain scarce data, or data from joint densities that do not have a corresponding optimal polynomial basis in the Askey scheme. As is remarked upon in Oladyshkin and Nowak, available distributions might be skewed or multi-modal \cite{OLADYSHKIN2012179}. To remedy this challenge, arbitrary Polynomial Chaos (aPC) has been proposed, in which the statistical moments of the joint density are used to derive the corresponding optimal basis \cite{AHLFELD20161}. 

{In this paper we marry aPC with recently proposed schemes that use sensitivity information to enrich the LSQ formulation that estimates the linear combination terms {\cite{se-pce, cfd_adjoint}}}. This sensitivity information can be harvested from the adjoint system at comparable cost to a single evaluation of $\mathcal{M}$, providing $n_u$ additional equations, besides the solution of $\mathcal{M}(\boldsymbol{\xi}^{(i)})$ at the $i$\textsuperscript{th} sample point in the uncertain space. In so doing, the computational cost of employing aPC in high dimensional problems is greatly reduced, scaling asymptotically as $n_u^{p-1}$. Using sensitivity information to enrich a surrogate model is an established strategy in UQ, with examples in the literature involving kriging \cite{han2013improving}, co-kriging \cite{chung2002using}, robust design \cite{iPCE_Barca}, and polynomial regression \cite{roderick2010polynomial} techniques. There are also instances of sensitivity/gradient information being used for PCE surrogate models, for instance, gradient information has been utilised as part of an $l1$ minimization process to identify the coefficients of a PCE \cite{guo2018gradient,jakeman2015enhancing,peng2016polynomial,Lockwood_Mavriplis_2013,Luchini2014AdjointAnalysis,RoderickFischer,Isukapalli_et_al_2000,PENG2016440,Guo_et_al_2018b}, among others. In what follows we outline the SeAr-PC method, and demonstrate its application to a set of high-dimensional synthetic test cases.

\section{Sensitivity enhancement for arbitrary Polynomial Chaos}
This section outlines the method for sensitivity enhanced arbitrary Polynomial Chaos (SeAr PC). The first subsection describes how the statistical moments of $f(\boldsymbol{\xi})$ are used to determine an optimal polynomial basis, while the second explains how sensitivity information can be used to enrich the LSQ formulation that estimates the PCE coefficients and consequently, greatly reduce the number of required function evaluations. Finally, the third section outlines the D-optimal procedure used to determine the Design of Experiments (DoE).

\subsection{Arbitrary Polynomial Chaos}
PCEs provide a spectral representation of an uncertain QoI, denoted as: 

\begin{align}
    \mathcal{M}(\boldsymbol{\xi})=\sum^{P}_{k=0}{\boldsymbol{\lambda}}_k\Psi_k(\boldsymbol{\xi}),
    \label{pce}
\end{align}
where $\boldsymbol{\lambda}\subseteq\Re^{P+1}$ is a set of deterministic coefficients and $\Psi$ a multivariate orthogonal polynomial, itself a product of a set of uni-variate, orthogonal polynomials, $\psi$, with:
\begin{align}
    \Psi_k(\boldsymbol{\xi})=\prod_{i=1}^{n_u}\psi_{I_{k,i}}(\boldsymbol{\xi}_i), \;\; k=0,\dots, P.
\end{align}
$I$ is a $(P+1) \times n_u$ index matrix, with the rows denoting the corresponding orders of the uni-variate polynomials for each term in the expansion. Xiu and Karniadakis \cite{xiu} proved that members of the Askey scheme of orthonormal polynomials are optimal for certain types of stochastic processes  such  as  the  uniform, normal, and exponential processes. However, should $f(\boldsymbol{\xi})$ not belong to this set of processes, then aPC can be used to determine the optimal uni-variate polynomials. Note that is is assumed within the PCE formulation that the elements of $\boldsymbol{\xi}$ are independent, i.e. $f(\boldsymbol{\xi}$)=$\prod_{i=1}^{n_u}f_i(\boldsymbol{\xi}_i)$, where $f_i$ is the Probability Density Function (PDF) for the $i$\textsuperscript{th} uncertain input parameter. 

A feature of PCEs is that the statistical moments of $f(\mathcal{M})$ can be expressed analytically using $\hat{\boldsymbol{\lambda}}$, with the mean and variance expressed as:
\begin{align}
    &\mathbb{E}(\mathcal{M})=\hat{\boldsymbol{\lambda}}_0, \label{eq:pce_mom}\\
    &\text{Var}(\mathcal{M})=\sum_{j=1}^P\gamma_j\hat{\boldsymbol{\lambda}}_, \nonumber
\end{align}
where $\gamma_j$ is a normalising constant estimated as:

\begin{align}
    \gamma_j=\mathbb{E}\big(\Psi_j^2(\xi)\big).
\end{align}

In this paper we calculate the optimal uni-variate polynomials, $\psi$, from the statistical moments of the input data. The remainder of this section summarises the main steps of the algorithm. The interested reader is referred to Ahfeld et al. \cite{AHLFELD20161} for a fuller description of the mathematical details of moment-based arbitrary Polynomial Chaos. The statistical moments of the input data can be provided directly to the algorithm if analytic expressions are available; or as a set of $n_s$ mono-dimensional samples $\{{\xi}^{(1)},\dots, {\xi}^{(n_s)}\}$, in which case the $i$\textsuperscript{th} raw moment, $\mu_i$, may be calculated as:

\begin{align}
    \mu_i=\sum_{l=1}^{n_s} (\xi^{(l)})^i.
\end{align}
The use of the statistical moments gives the approach flexibility and allows sparse input data to be efficiently handled. Having calculated the first $2p+1$ statistical moments, the Hankel matrix of moments is formed:

\begin{align}
    M=\begin{bmatrix}
        \mu_0 & \mu_1 & \ldots& \mu_p\\
        \mu_1 & \mu_2 & & \mu_{p+1}\\
        \vdots &  & \ddots& \\
        \mu_p & \mu_{p+1} & ...& \mu_{2p}
    \end{bmatrix}.
    \label{samba1}
\end{align}
The Hankel matrix is positive definite, therefore the Cholesky decomposition may be computed, yielding the upper triangular matrix $R$, where $M=R^\top R$. The Mysovskih Theorem states that the entries of the inverse matrix $R^{-1}$ form an orthogonal system of polynomials:

\begin{align}
    \psi_j=\sum_{k=0}^{p}\underbar{r}_{k+1,j+1}\,\xi^k,
    \label{samba2}
\end{align}
where $\underbar{r}_{i,j}$ denotes the elements of the (upper triangular) inverse matrix, $R^{-1}$.  This inverse can be computed directly or alternatively, through the analytic formulas of Rutishauser, the polynomial basis can be derived from a set of three term recurrence relations.

\subsection{Sensitivity Enhanced arbitrary Polynomial Chaos}
Having determined the optimal set of orthogonal polynomials for each uncertain input, a least squares minimisation is formulated to estimate the PCE coefficients, $\boldsymbol{\lambda}=[\lambda_0,..., \lambda_P]^\top$:

\begin{equation}
\hat{\boldsymbol{\lambda}}=\min_{\boldsymbol{\lambda}}\; ( {Q} - \boldsymbol{\psi} \boldsymbol{\lambda})^\top\,W\, ( {Q} - \boldsymbol{\psi} \boldsymbol{\lambda} ),
\label{min}
\end{equation}
where $Q=[\mathcal{M}(\boldsymbol{\xi}^{(1)}),\dots, \mathcal{M}(\boldsymbol{\xi}^{(q)})]^\top$ collects the evaluations of $\mathcal {M}$ for a set of $q$ samples $\{\boldsymbol{\xi}^{(1)},\dots, \boldsymbol{\xi}^{(q)}\}$. $W$ is a diagonal matrix of weights that scale the relative contribution of each sample. The choice of this weighting matrix is discussed further in the next subsection. The $q \times (P+1)$ measurement matrix, $\boldsymbol{\psi}$, collects the evaluations of each term of the polynomial basis at the sample points:

\begin{equation}
\boldsymbol{\psi} = 
\begin{bmatrix}
\Psi_0 \left(\boldsymbol{\xi}^{(1)}\right)   &  \dots  & \Psi_P\left(\boldsymbol{\xi}^{(1)}\right)  \\
  \vdots     &  \ddots & \vdots    \\
\Psi_0\left(\boldsymbol{\xi}^{(q)}\right)  &  \dots   & \Psi_P\left(\boldsymbol{\xi}^{(q)}\right)   \\
\end{bmatrix}.
\label{eq:psi_matrix}
\end{equation}
Solutions to \eqref{min}, $\hat{\boldsymbol{\lambda}}$, satisfy the normal equation:

\begin{align}
    \left( \boldsymbol{\psi}^{\top} {W} \boldsymbol{\psi} \right) \hat{\boldsymbol{\lambda}} =   \boldsymbol{\psi}^{\top} {W} {Q}.
\end{align}
In order for \eqref{min} to be well conditioned, at least $P+1$ samples are required but preferably $q\gg P+1$. Given that each additional sample requires an evaluation of $\mathcal{M}$, which could be an expensive computer code, there is great advantage in augmenting the system with additional information that does not require additional model evaluations. In SeAr PC, this is done by augmenting \eqref{min} with the sensitivities of $\mathcal{M}(\boldsymbol{\xi})$ to each element of $\boldsymbol{\xi}$. {This section summarises the main steps through which sensitivity information can be incorporated in the PCE formulation. For further details on sensitivity enhancement for Askey scheme PCE bases, the reader is referred to Kantarakias and Papadakis {\cite{se-pce}}}.

\begin{figure}
\begin{center}
\includegraphics[width=0.6\textwidth]{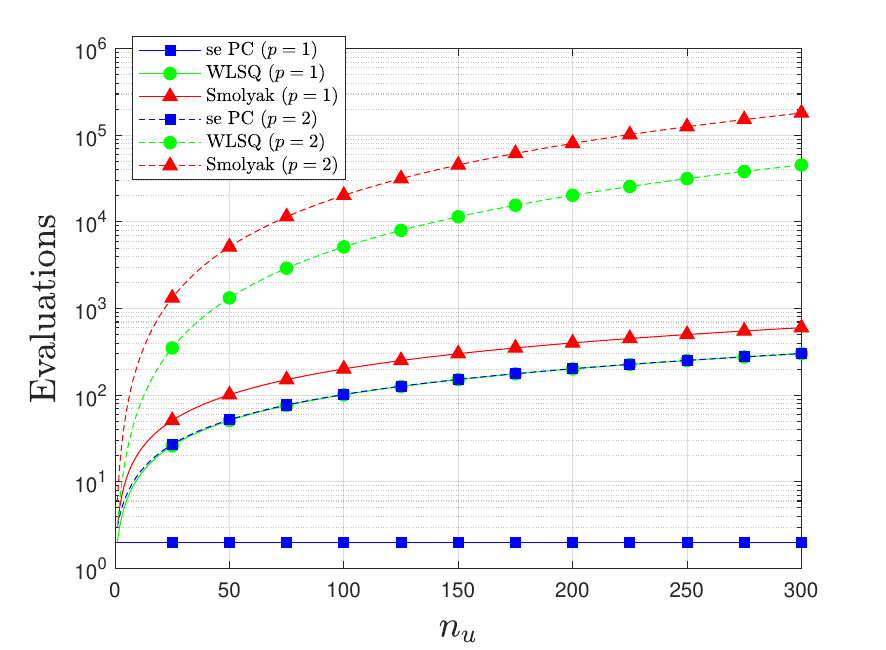}
\end{center}
\caption{Minimum number of model evaluations required as a function of $n_u$ and $p$ for sensitivity enhanced PC, weighted least squares, and Smolyak sampling grid.}
\label{fig:smolyak} 
\end{figure}

Including sensitivity information enables $n_u$ additional equations may be generated for each sample. Taking the derivative of \eqref{pce} with respect to the $i$\textsuperscript{th} element of $\boldsymbol{\xi}$ yields: 

\begin{equation}
\frac {d \mathcal{M}}{d \boldsymbol{\xi}_i^{(j)}}  =\sum_{k=0}^{P} \lambda_{k} \frac{\partial \Psi_k}{ \partial \boldsymbol{\xi}_i^{(j)}}, \quad j=1,\dots,q.
\label{eq:sepce01}
\end{equation}
Each component therefore yields the block of equations:
\begin{equation}
\frac{d {{Q}}}{d \boldsymbol{\xi}_i}=
\frac{\partial \boldsymbol{\psi} }{\partial \boldsymbol{\xi}_i}
\boldsymbol{\lambda}, \quad i=1,\dots,n_u,
\label{eq:sepce02}
\end{equation}
with 
\begin{align}
 \frac{d {{Q}}}{d \boldsymbol{\xi}_i} =\left [\frac{d \mathcal{M}}{d \boldsymbol{\xi}_i^{(1)} }, \dots, \frac{d \mathcal{M}}{d \boldsymbol{\xi}_i^{(q)} } \right]^\top,
 \label{dq}
 \end{align} 
and 
\begin{equation}
\frac{\partial \boldsymbol{\psi} }{\partial \boldsymbol{\xi}_i}=
\begin{bmatrix}
\frac{\partial \Psi_0 \left(\boldsymbol{\xi}^{(1)}\right) }{\partial \boldsymbol{\xi}_i^{(1)}}   &  \dots  & \frac{\partial \Psi_P \left(\boldsymbol{\xi}^{(1)}\right)}{\partial \boldsymbol{\xi}_i^{(1)}}   \\
  \vdots     &  \ddots & \vdots    \\
\frac{\partial \Psi_0 \left(\boldsymbol{\xi}^{(q)} \right)}{\partial \boldsymbol{\xi}_i^{(q)}}   &  \dots   & \frac{\partial \Psi_P \left(\boldsymbol{\xi}^{(q)} \right)}{\partial \boldsymbol{\xi}_i^{(q)}}
\end{bmatrix}.\\
\label{eq:nabla_psi_matrix}
\end{equation}
Elements of the matrix in \eqref{eq:nabla_psi_matrix} can be derived analytically from the expressions for the optimal polynomials. On the other hand, the sensitivities in \eqref{dq} can be found through the adjoint formulation of the system. Repeating this process for each uncertain input yields the $(n_u+1)q \times 1$ block column vector that collects the model evaluations and sensitivity information:
\begin{equation}
\begin{aligned}
& {G}  =  \left [{Q}, \frac{d {{Q}}}{d \boldsymbol{\xi}_1}, \dots, \frac{d {{Q}}}{d \boldsymbol{\xi}_{n_u}} \right ]^{\top}= \\
 & \left[\mathcal{M} \left(\boldsymbol{\xi}^{(1)}\right),\dots,\mathcal{M} \left(\boldsymbol{\xi}^{(q)}\right), \frac{d \mathcal{M}}{d \boldsymbol{\xi}_1^{(1)} }, \dots, \frac{d \mathcal{M}}{d \boldsymbol{\xi}_1^{(q)}}, \dots, 
\frac{d \mathcal{M}}{d \boldsymbol{\xi}_{n_u}^{(1)}},\dots,\frac{d \mathcal{M}}{d \boldsymbol{\xi}_{n_u}^{(q)}}  \right]^\top 
\end{aligned}
\label{eq:expanded_G}
\end{equation}
Similarly a $(n_u+1)q \times (P+1)$ row block matrix can be defined for the multivariate polynomials and their derivatives:
\begin{equation}
\boldsymbol{\phi}=
\begin{bmatrix}
\boldsymbol{\psi} \\  \frac{\partial \boldsymbol{\psi} }{\partial \boldsymbol{\xi}_1} \\ \vdots \\ \frac{\partial \boldsymbol{\psi} }{\partial \boldsymbol{\xi}_{n_u}} 
\end{bmatrix}=
 \begin{bmatrix}
\Psi_0 \left(\boldsymbol{\xi}^{(1)}\right)   &  \dots  & \Psi_P\left(\boldsymbol{\xi}^{(1)}\right)  \\
  \vdots     &  \ddots & \vdots    \\
\Psi_0\left(\boldsymbol{\xi}^{(q)}\right)  &  \dots   & \Psi_P\left(\boldsymbol{\xi}^{(q)}\right)   \\
\frac{\partial \Psi_0 \left(\boldsymbol{\xi}^{(1)}\right) }{\partial \boldsymbol{\xi}_1^{(1)}}   &  \dots  & \frac{\partial \Psi_P \left(\boldsymbol{\xi}^{(1)}\right)}{\partial \boldsymbol{\xi}_1^{(1)}}   \\
  \vdots     &  \ddots & \vdots    \\
\frac{\partial \Psi_0 \left(\boldsymbol{\xi}^{(q)} \right)}{\partial \boldsymbol{\xi}_1^{(q)}}   &  \dots   & \frac{\partial \Psi_P \left(\boldsymbol{\xi}^{(q)} \right)}{\partial \boldsymbol{\xi}_1^{(q)}} \\
  \vdots     &  \ddots & \vdots    \\
\frac{\partial \Psi_0 \left(\boldsymbol{\xi}^{(1)} \right)}{\partial \boldsymbol{\xi}_{n_u}^{(1)}}   &  \dots   & \frac{\partial \Psi_P \left(\boldsymbol{\xi}^{(1)} \right)}{\partial \boldsymbol{\xi}_{n_u}^{(1)}}   \\
 \vdots     &  \ddots & \vdots \\
\frac{\partial \Psi_0 \left(\boldsymbol{\xi}^{(q)} \right)}{\partial \boldsymbol{\xi}_{n_u}^{(q)}}   &  \dots   & \frac{\partial \Psi_P \left(\boldsymbol{\xi}^{(q)} \right)}{\partial \boldsymbol{\xi}_{n_u}^{(q)}}   
\end{bmatrix}
 \label{eq:expanded_phi}
\end{equation}
Finally, these block matrices can be used to redefine the least squares optimisation for $\hat{\boldsymbol{\lambda}}$: 
\begin{equation}
\hat{\boldsymbol{\lambda}}=\min_{\boldsymbol{\lambda}}  \left({G}- \boldsymbol{\phi} \boldsymbol{\lambda}\right)^\top  {W}' \left({G}- \boldsymbol{\phi} \boldsymbol{\lambda}\right),
\label{eq:sepce03}
\end{equation}
where $ {W'}$ is a block diagonal weighting matrix, consisting of $n_u+1$ blocks of ${W}$. Figure \ref{fig:smolyak} illustrates the benefits of this formulation in terms of the number of required samples. A necessary condition for $\hat{\boldsymbol{\lambda}}$ to be a unique solution to the least squares problem is that the measurement matrix must have at least as many rows as required coefficients ($P+1$). In consequence, $P+1$ model evaluations are necessary for the standard formulation of a PCE, while $(P+1)/(n_u+1)$ are required for a PCE with sensitivity enhancement. This dependence is plotted in Figure \ref{fig:smolyak} for $p=1$ and $p=2$ for an oversampling ratio of 2. Note that in the case of $p=1$, a consequence of the sensitivity enhancement is that the number of required samples is independent of $n_u$. For comparison, the necessary number of model evaluations of weighted least squares and a Smolyak sampling grid are plotted.  

\subsection{D-optimal Design of Experiments}
A least squares solution to \eqref{eq:sepce03} is sought to estimate $\hat{\boldsymbol{\lambda}}$. The coherence parameter, $\mu_c$, is understood to have a significant effect on the convergence and stability of this solution (see, e.g. \cite{cohen2013stability}). The sub-script ``c" is used to distinguish this quantity from the statistical moments. $\mathcal{M}(\boldsymbol{\xi})$ is evaluated at $q$ locations in $\Xi$. It is therefore desirable to identify $q$ samples that are coherence optimal, in the sense that they minimise $\mu_c$, and use this as the Design of Experiments (DoE) for the uncertainty analysis. In Hampton and Doostan \cite{doostan2}, $\mu_c$ is defined for a PCE, using the general result of Cohen et al \cite{cohen2013stability}, as:
\begin{align}
    \mu_c=\text{sup}\sum_{j=1}^P|\mathcal{W}(\boldsymbol{\xi})\Psi_j(\boldsymbol{\xi})|^2, 
\end{align}
where $\mathcal{W}(\cdot)$ is a weighting function for $f(\boldsymbol{\xi})$. $\mathcal{W}(\cdot)$ has been derived for certain polynomial bases, for instance for (physicists') Hermite polynomials the function \cite{doostan2}:
\begin{align}
     \mathcal{W}(\boldsymbol{\xi})=\text{exp}\bigg(-\frac{1}{4}\|\boldsymbol{\xi}\|^2\bigg),
\end{align}
is used, while:
\begin{align}
     \mathcal{W}(\boldsymbol{\xi})=\prod^{n_u}_{i=1}(1-\boldsymbol{\xi}_i^2)^{1/4},
\end{align}
can be used for Legendre polynomials. However, this paper is concerned with broadening the applicability of polynomial chaos to stochastic processes without a corresponding optimal polynomial in the Askey scheme. To this end, we suggest use the function proposed in  Hampton and Doostan \cite{doostan1}:

\begin{align}
    \mathcal{W}(\boldsymbol{\xi})=c^{-1}\mathcal{B}(\boldsymbol{\xi})^{-1},
\end{align}
where $c$ is a normalising constant and $\mathcal{B}(\cdot)$ is itself a function of $\boldsymbol{\xi}$ and the multivariate polynomial basis from aPC:
\begin{align}
    \mathcal{B}(\boldsymbol{\xi})=\sqrt{\sum_{j=1}^P|\Psi_j(\boldsymbol{\xi})|^2}.
\end{align}
$\mathcal{B}^2(\boldsymbol{\xi})$ represents a uniformly least upper bound on the sum of squares of the multi-variate polynomial basis. The normalising constant, $c$, is defined as:
\begin{align}
    c^2=\int_{\Xi}f(\boldsymbol{\xi})\mathcal{B}^2(\boldsymbol{\xi})d\boldsymbol{\xi}. 
\end{align}
However, as with \eqref{mean}, there is no analytic expression for this integral. There are therefore two alternative coherent sampling strategies: one is to use Markov Chain Monte Carlo to sample from the coherent optimal probability density:
\begin{align}
    f_\mathcal{M}(\boldsymbol{\xi})=c^2\mathcal{B}^2(\boldsymbol{\xi})f(\boldsymbol{\xi}). 
\end{align}
Alternatively, a pool of samples can be drawn from $f(\boldsymbol{\xi})$ using Monte Carlo sampling and from this pool, a coherence optimal set of $q$ samples can be selected. In this paper the former option is chosen, with the pivoted QR decomposition used to identify the coherence optimal points. 

\subsubsection{Coherence optimal sampling through the pivoted QR decomposition}
In this sub-section the greedy algorithm used to determine the DoE is summarised. {Following the approach of Kantarakias and Papadakis {\cite{se-pce}}, the DoE for a sensitivity enhanced PCE is determined through an optimisation process}. The pivoted QR decomposition algorithm \cite{qr1,qr2} is used to determine an optimal set of {$q$} locations in {$\Xi$} at which to evaluate {$\mathcal{M}$}. These are selected from a pool of $n_s$ samples of the joint density, with $n_s\gg q$. {Key relations of the method used to determine the optimal DoE are relayed here. We refer the reader to Kantarakias and Papadakis {\cite{se-pce}} for the original derivations.}

In SeAr PC a two step process is followed: first the $q_a$ coherence optimal sampling locations for aPC without the sensitivity enhancement are identified. For an oversampling ratio of $n_o$, this corresponds to $q_a=n_o(P+1)$ samples. From these samples, the subset of the $q$ most important are chosen for the DoE, where:

\begin{align}
    q=\frac{n_o(P+1)}{n_u+1}=\frac{q_a}{n_u+1}. 
\end{align}

Restating the least squares optimisation for aPC without sensitivity enhancement, \eqref{min}, as a matrix equation yields: 

\begin{align}
   W^{\frac{1}{2}} Q=W^{\frac{1}{2}}\boldsymbol{\psi}{\boldsymbol{\lambda}}+\boldsymbol{\epsilon}, 
   \label{interp2}
\end{align}
where $\boldsymbol{\epsilon}\in\Re^{P+1}$ represents the interpolation error, assumed to be independent Gaussian noise with standard deviation $\eta$, i.e. $\boldsymbol{\epsilon}_i\approx N({0},\eta)$. The diagonal matrix, $W$, weights the contribution of each sample a sample in $\Xi$. The entries in $W$ correspond to an evaluation of $\mathcal{W}(\cdot)$ for each of the $n_s$ samples, i.e.:

\begin{align}
    W_{ii}=\mathcal{W}(\xi^{(i)}), \; i=1,\dots, n_s.
    \label{eq:w}
\end{align}
The error covariance matrix is used as a measure of the coherence, which given \eqref{interp2}, can be shown to be \cite{se-pce}:
\begin{align}
    \text{Var}(\boldsymbol{\lambda}-\hat{\boldsymbol{\lambda}})=\eta^2(\boldsymbol{\psi}^\top W \boldsymbol{\psi})^{-1}, 
    \label{eq:det}
\end{align}
Minimising this quantity is equivalent to maximising the determinant of the RHS. What is desired is the row selection matrix, which we denote $P_a$, that selects the rows of the measurement matrix that correspond to coherence optimal samples i.e. 

\begin{align}
    \boldsymbol{\psi}_a= P_a\boldsymbol{\psi},
\end{align}

where at each row of $P_a$ all elements are 0, except the element at the column that corresponds to the selected sampling point, which takes the value of 1. The maximisation problem then becomes:

\begin{align}
    P_a&=\underset{P_a}{\text{argmax}}\;\text{det}\big[(P_a W^{\frac{1}{2}} \boldsymbol{\psi})^\top (P_a W^{\frac{1}{2}} \boldsymbol{\psi})\big]. \label{pivotQR}\\
    &=\underset{P_a}{\text{argmax}}\;\text{det}\big[P_a W^{\frac{1}{2}} \boldsymbol{\psi}\big] \nonumber
\end{align}
This optimisation problem is solved using the pivoted QR decomposition: 
\begin{align}
    (W^{\frac{1}{2}}\boldsymbol{\psi})^\top P_a^\top=Q_aR_a, 
\end{align}
where $R_a\in\Re^{q_a\times q_a}$ is a upper diagonal matrix. The subscripts ``$a$" are used to differentiate these matrices from the definitions in previous sections. $P_a$ is selected such that the diagonal elements of $R_a$ are ordered in descending order, i.e.: 
\begin{align}
    |R_{a_{11}}|\geq |R_{a_{22}}|\geq \dots \geq |R_{a_{q_a,q_a}}|.
\end{align}
The absolute value of the determinant is given by the product of the diagonal entries of $R_a$:
\begin{align}
    \big|\text{det}\big(W^{\frac{1}{2}}\boldsymbol{\psi}\big)\big|=\prod_{m=1}^{q}|R_{mm}|.
\end{align}

\begin{figure}
\begin{center}
\includegraphics[width=0.49\textwidth]{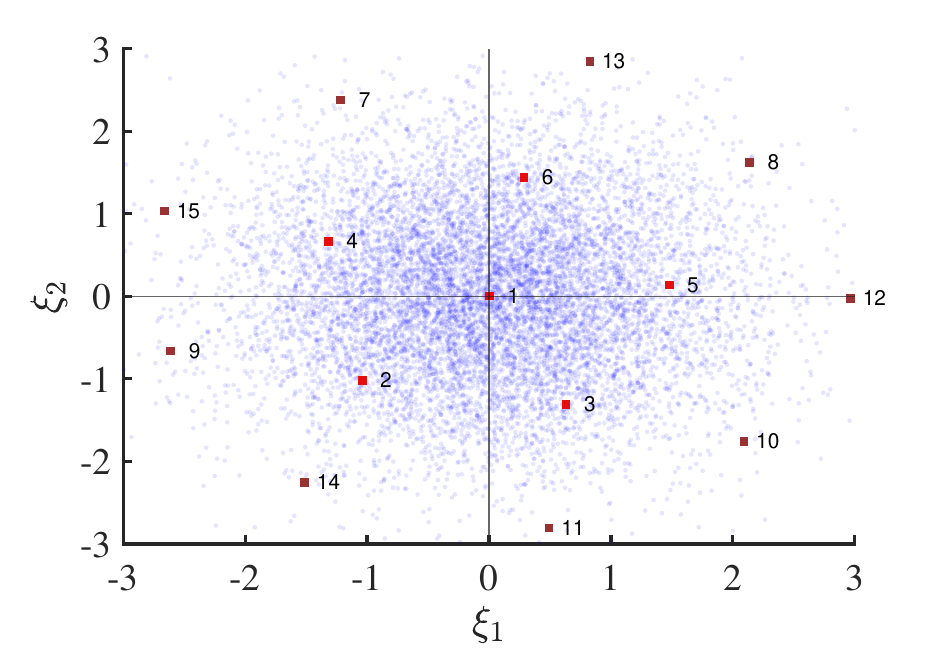}
\includegraphics[width=0.49\textwidth]{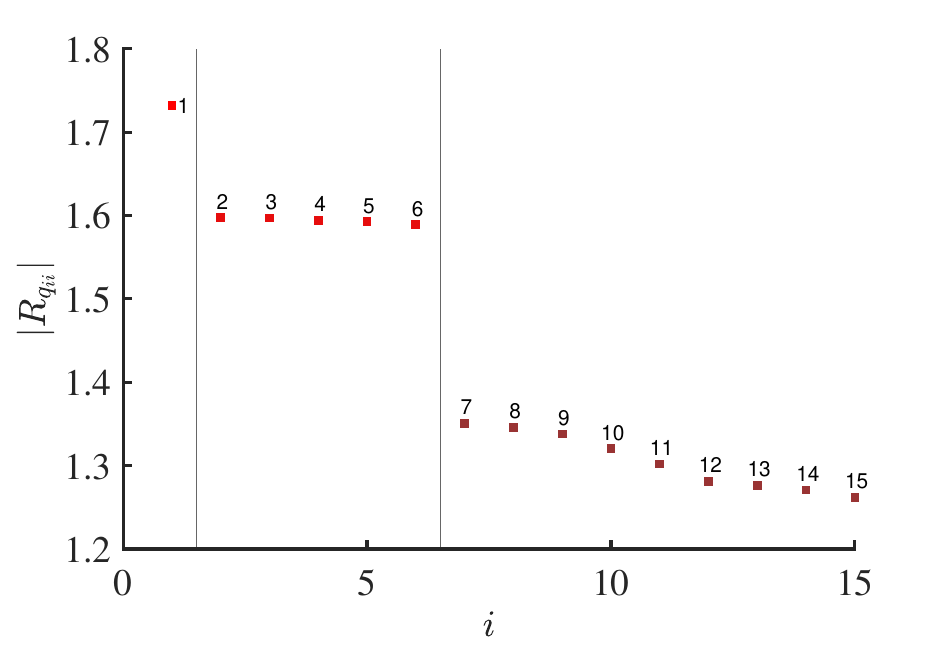}
\includegraphics[width=0.49\textwidth]{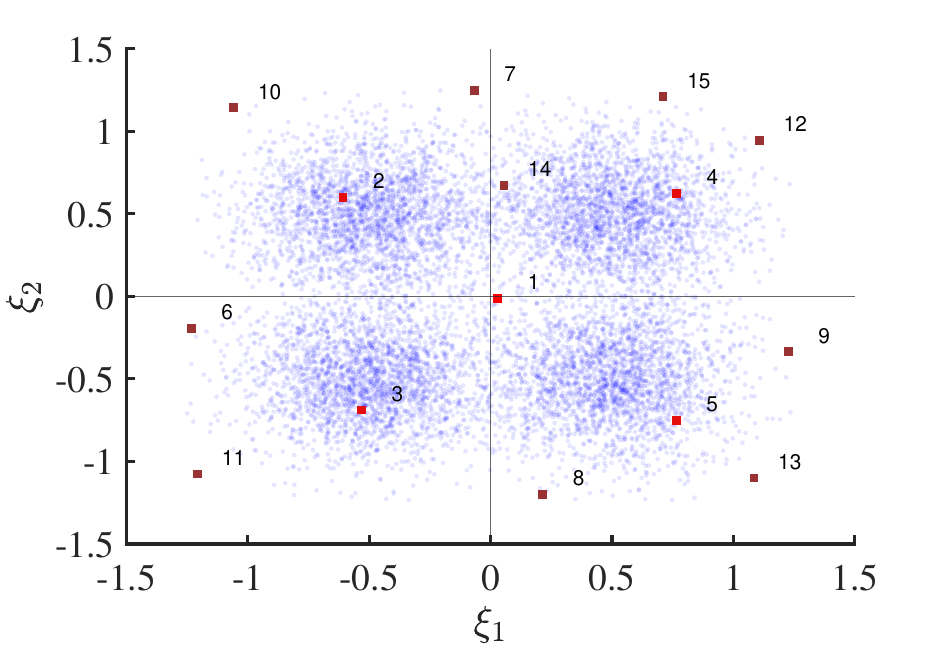}
\includegraphics[width=0.49\textwidth]{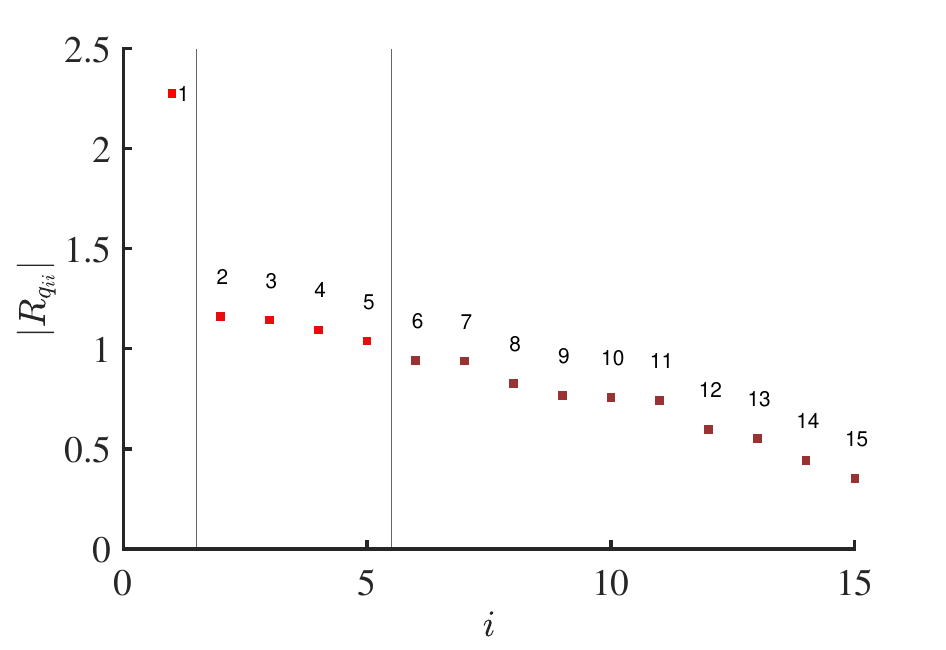}
\end{center}
\caption{Illustration of the D-optimal sampling points and corresponding diagonal entries in $|R_a|$ for a set of independent Gaussian distributions (top) and Gaussian mixture (bottom).}
\label{fig:samp_grid} 
\end{figure}

Having identified $q_a$ optimal sampling points in $\Xi$ for aPC without the sensitivity enhancement, the first $q$ points are selected for SeAr-PC, recalling that fewer points are required if sensitivity information is included. Figure \ref{fig:samp_grid} demonstrates the D-optimal sampling method for two $n_u=2$ uncertain spaces, with $n_o=1$ and $p=4$. The top left panel illustrates Monte Carlo samples from a bi-variate Gaussian distribution, with the standard normal distribution, $N(0,1)$, in each dimension used as the joint density. The locations of the D-optimal samples (with $p=4$) are indicated in red. The right panel of Figure \ref{fig:samp_grid} plots these sampling points against the corresponding diagonal entry of $R_a$. The sampling points are loosely organised in concentric rings around the origin, reflected by the clustering of values in the top right panel. The bottom two panels display the same plots for a joint density consisting of two independent Gaussian mixtures $N(-0.5,0.25)$ and $N(0.5,25)$. 

The SeAr PC algorithm may be summarised as follows:

\vspace{0.5cm}
\noindent {\bf Algorithm 1 (Sensitivity Enhanced Arbitrary Polynomial Chaos)}: \\
Inputs:  Joint density $f(\boldsymbol{\xi})$, computational model $\mathcal{M}(\boldsymbol{\xi})$\\
Outputs: Probability density $f(\mathcal{M}(\boldsymbol{\xi}))$
\begin{enumerate}

\item Draw $n_s$ Monte Carlo samples from $f(\boldsymbol{\xi})$

\item Compute statistical moments of the samples

\item Determine optimal multi-variate polynomial basis, $\{\Psi_k\}_{k=0}^P$, \eqref{samba1} \& \eqref{samba2}

\item Find $q_a$ optimal sampling locations through pivoted QR decomposition \eqref{pivotQR}, from these select the $q$ most significant as the DoE

\item Evaluate $\mathcal{M}(\boldsymbol{\xi})$ at these locations 

\item Perform least squares fitting for the PCE weights, $\{\hat{\boldsymbol{\lambda}}_k\}_{k=0}^P$ through \eqref{eq:sepce03}

\item Evaluate PCE for the $n_s$ Monte Carlo samples, fit $f(\mathcal{M}(\boldsymbol{\xi}))$ with kernel density smoothing

\end{enumerate}
\vspace{0.1cm}

\section{Synthetic Test Cases}
Having outlined the main steps of SeAr PC, we demonstrate the application of the algorithm to a set of synthetic test cases featuring multi-modal, fat-tailed, and truncated probability distributions. These densities do not have a corresponding optimal polynomial in the Askey scheme. Additional complication is introduced in that the uncertain input spaces are high-dimensional, with $n_u\geq 10$ for the three cases. The presented algorithm is baselined against an approach that estimates $\hat{\boldsymbol{\lambda}}$ through weighted least squares (WLSQ), with an aPC multi-variate polynomial basis. 

For the first two tests, we perform UQ for two test functions: a smooth, cubic function, defined as:

\begin{align}
    \mathcal{M}(\boldsymbol{\xi})=1+\boldsymbol{\xi}_1+\frac{1}{n_u}\sum_{i=1}^{n_u} \boldsymbol{\xi}_i^3,
\end{align}
and a non-linear sinusoidal function, previously used in Ahlfeld et al \cite{AHLFELD20161}:
\begin{align}
    \mathcal{M}(\boldsymbol{\xi})=\sum_{i=1}^{n_u}\textrm{sin}(\boldsymbol{\xi}_i-0.5).
\end{align}
We compare the convergence of SeAr-PC against the WLSQ for two different probability distributions (sampled independently in $n_u=10$ dimensions) with no optimal uni-variate polynomials in the Askey scheme. 



\subsection{Uncertainty Quantification with Multi-modal Input Probability Distributions}
\begin{figure}
\begin{center}
\includegraphics[width=0.49\textwidth]{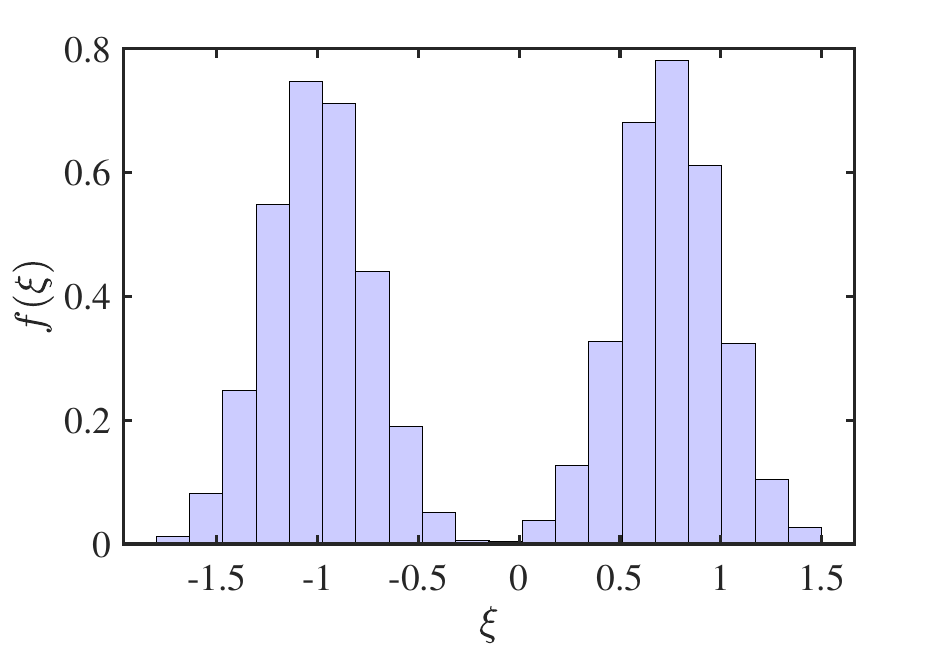}
\includegraphics[width=0.49\textwidth]{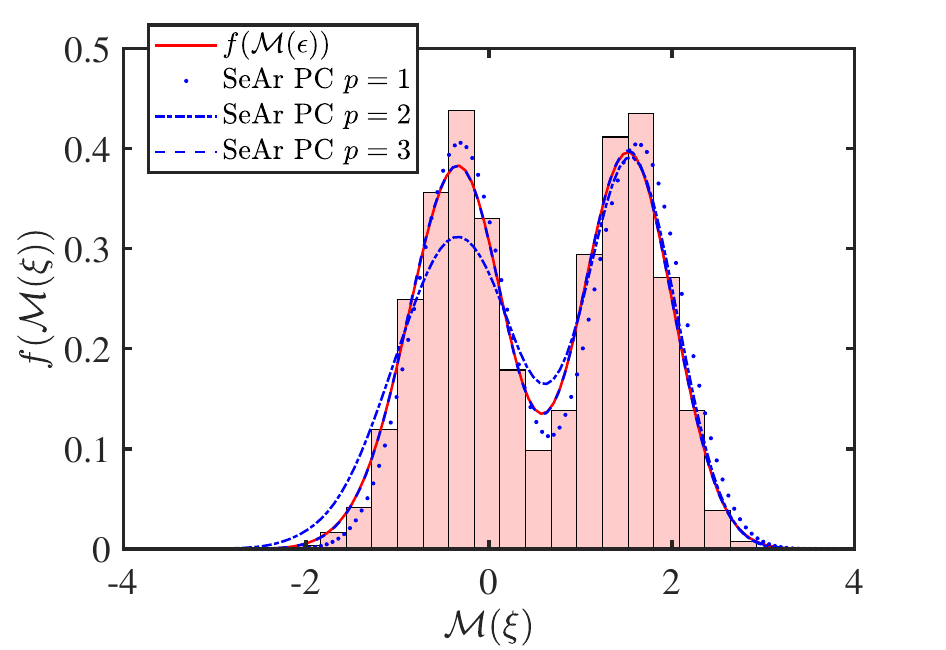}
\includegraphics[width=0.49\textwidth]{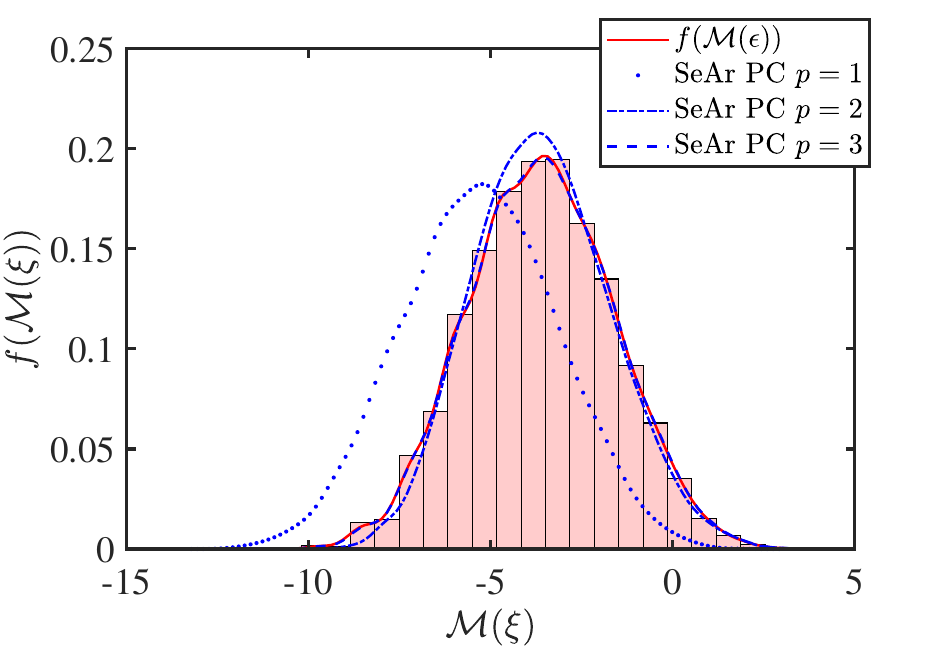}

\end{center}
\caption{Histogram of uni-variate Monte Carlo samples from the multi-modal joint density (top left) with the corresponding histogram found by Monte Carlo sampling for the cubic test function (top right) and non-linear function (bottom). Also indicated are the kernel density estimate of this distribution, together with those of SeAr PC for increasing polynomial order.}
\label{fig:test_inps} 
\end{figure}

Multi-modal probability distributions can frequently occur in physical systems and real-world datasets. Such distributions can represent a challenge for PCEs as they can require a high-order polynomial to represent them accurately using a non-optimal basis (see, e.g. \cite{nouy, pepper_mmodal}). 10,000 Monte Carlo samples were drawn from the mixture of the two Gaussians $N(-1.0, 0.25)$ and $N(0.75,0.25)$ and provided to the algorithm as histograms of uni-variate Monte Carlo samples. One of these histograms is displayed on the top left panel of Figure \ref{fig:test_inps}. The right panel displays the histogram for the QoI, found through Monte Carlo sampling. The kernel density estimate of the histogram is also displayed, together with the SeAr PC estimate as the order of the PCE, $p$, is increased. By $p=3$ the kernel density estimates of $f(\mathcal{M}(\boldsymbol{\xi}))$ matched exactly. For all experiments an oversampling ratio of 2 was chosen.%

\begin{figure}
\begin{center}
\includegraphics[width=0.49\textwidth]{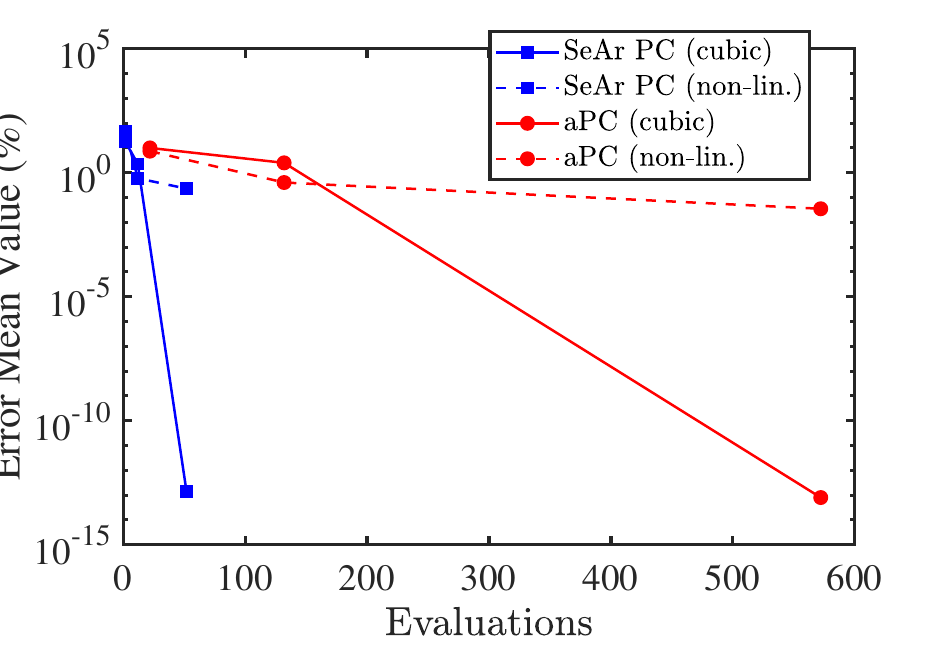}
\includegraphics[width=0.49\textwidth]{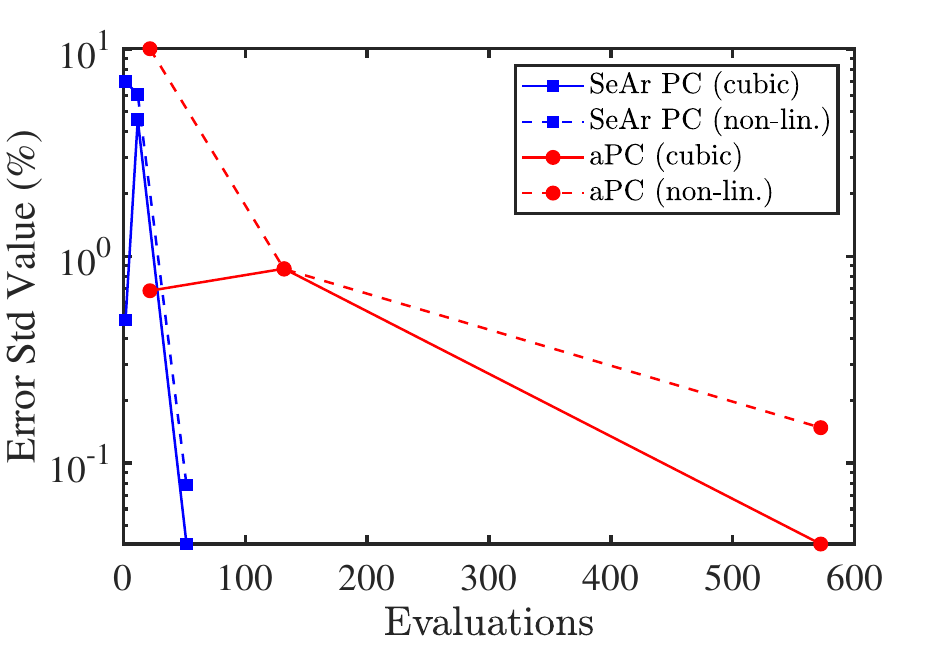}
\includegraphics[width=0.49\textwidth]{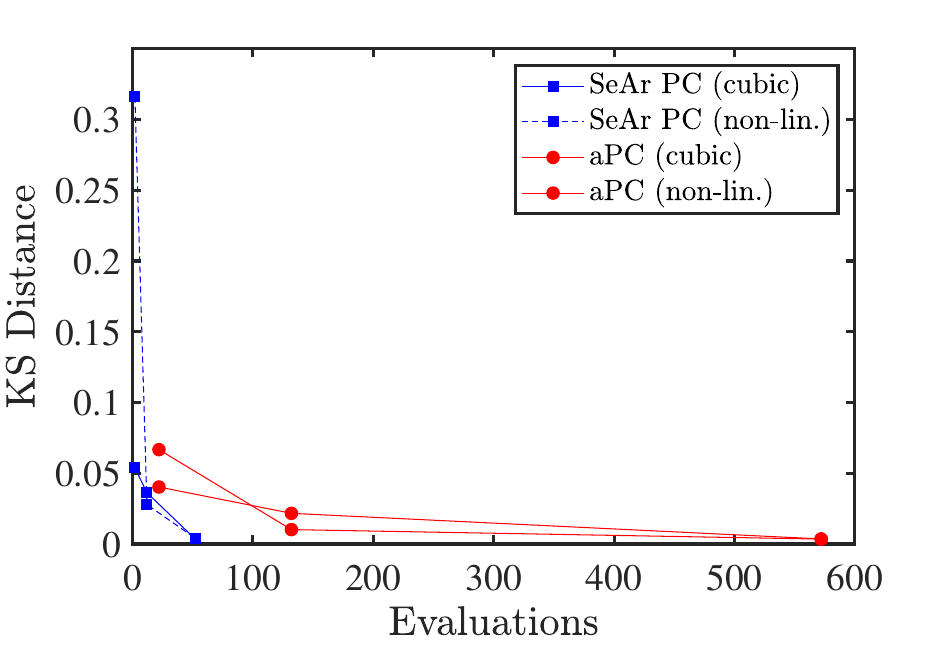}
\end{center}
\caption{Convergence of the first two statistical moments (top) and KS distance (bottom) with order/number of evaluations for the generalized extreme value test case.}
\label{fig:mmodal} 
\end{figure}

Figure \ref{fig:mmodal} shows the convergence of the first two statistical moments for each of the methods, plotted against number of evaluations, while the bottom panel indicates the convergence of the Kolmogorov-Smirnov (KS) distance. The solid lines indicate the application to the cubic function, the dotted lines the non-linear function. The square markers indicate the polynomial order. Data from Figure \ref{fig:mmodal} is tabulated in Table 1.  The plots demonstrate the benefits of the sensitivity enhancement in terms of the number of required samples, with the sensitivity enhanced implementations requiring many fewer samples to achieve the same accuracy as the corresponding WLSQ implementation. 


\begin{table}[!htbp]
\centering
\caption{Convergence of SeAr PC and aPC, for increasing polynomial order, for the generalized extreme value joint density. Results for the cubic and non-linear test functions are displayed. }
\begin{tabular}{*1c | *6c}
\toprule
Test case & Method & $p$ & $q$ & $\Delta \mu \,(\%)$ & $\Delta \sigma \,(\%)$ & KS dist. \\
\midrule
\multirow{6}{*}{Cubic function} &
\multirow{3}{*}{SeAr PC} &
1 &2 &17.22 &0.49 & 0.054\\
& & 2& 12& 2.29& 4.57& 0.036\\
& & 3& 52& 1.39e-13& 0.04& 0.004\\
& \multirow{3}{*}{aPC (WLSQ)}  &
 1 & 22& 9.98& 0.68& 0.040\\
& &2 &132 & 2.50 & 0.87&0.022 \\
&&3 &572 & 7.95e-14 &0.04 &0.004 \\
\midrule
\multirow{6}{*}{Non-lin. function} &
\multirow{3}{*}{SeAr PC} &
1 &2 & 46.71  & 6.98 & 0.316 \\
& & 2&12 & 0.60 & 6.05 & 0.028\\
& & 3& 52& 0.23 & 0.08  & 0.004\\
& \multirow{3}{*}{aPC (WLSQ)}  &
 1 &22 & 7.39 & 10.07 & 0.067\\
& &2 &132& 0.40 & 0.87 & 0.010\\
&&3 &572 & 0.04 & 0.15 & 0.004\\
\bottomrule
\end{tabular}
\end{table}

\subsection{Fat tailed, non-Askey scheme inputs}
\begin{figure}
\begin{center}
\includegraphics[width=0.49\textwidth]{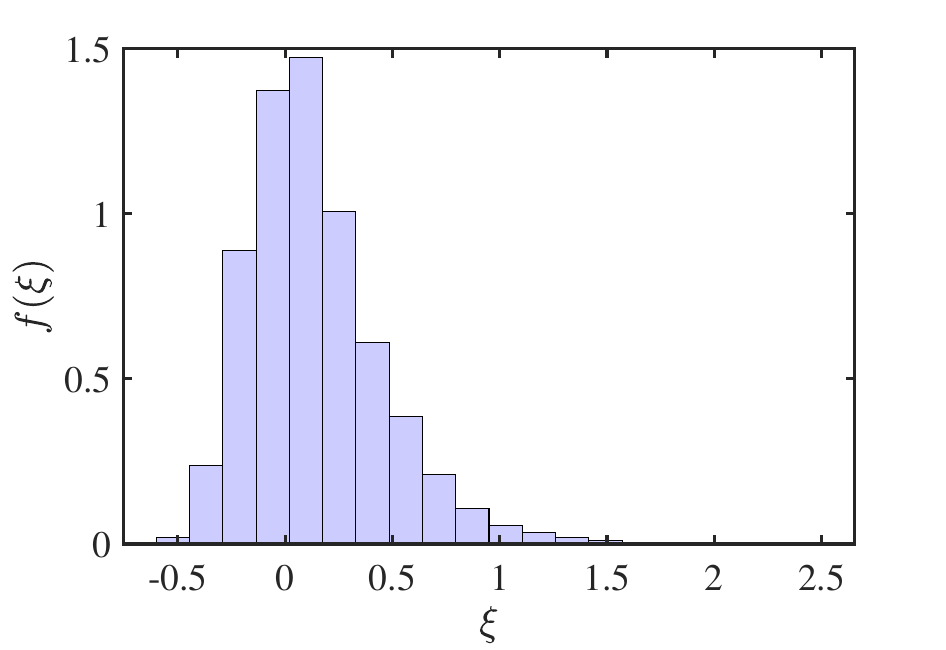}
\includegraphics[width=0.49\textwidth]{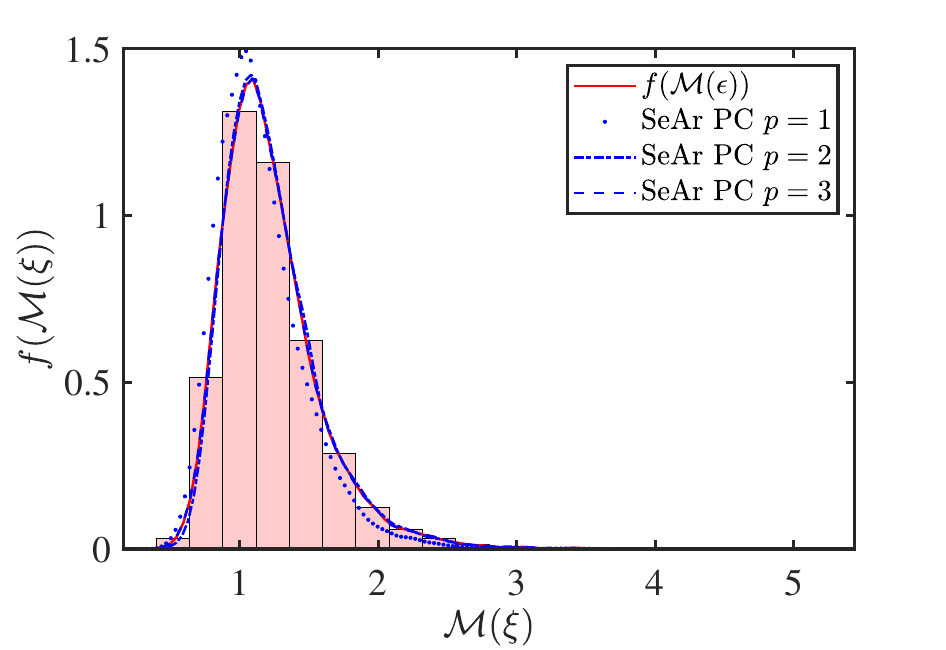}
\includegraphics[width=0.49\textwidth]{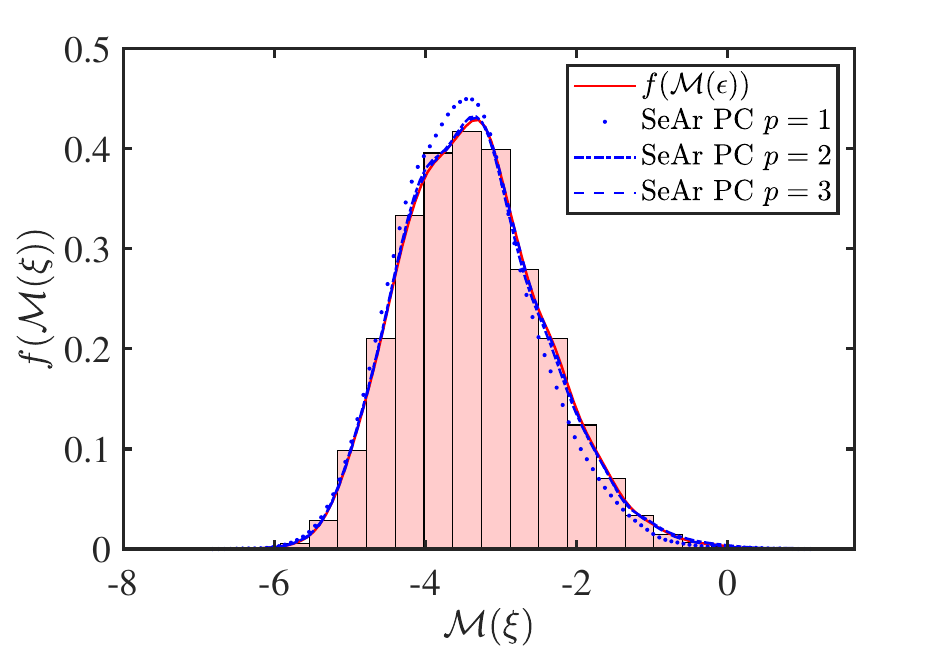}

\end{center}
\caption{Histogram of uni-variate Monte Carlo samples from the generalized extreme value joint density (top left), with the corresponding histogram found by Monte Carlo sampling for the cubic test function (top right) and non-linear function (bottom). Also indicated are the kernel density estimate of this distribution, together with those of SeAr PC for increasing polynomial order.}
\label{fig:test_inps_gev} 
\end{figure}

The analysis of fat-tailed stochastic processes is of great importance in reliability analysis. In this case $n_u$ independent Generalised Extreme Value distributions were chosen, with parameters $GEV(0,0.25,0)$. This is a useful case from an industrial perspective due to the heavy right tail of the distribution (visualised in the top left panel of Figure \ref{fig:test_inps_gev}). 

\begin{figure}
\begin{center}
\includegraphics[width=0.49\textwidth]{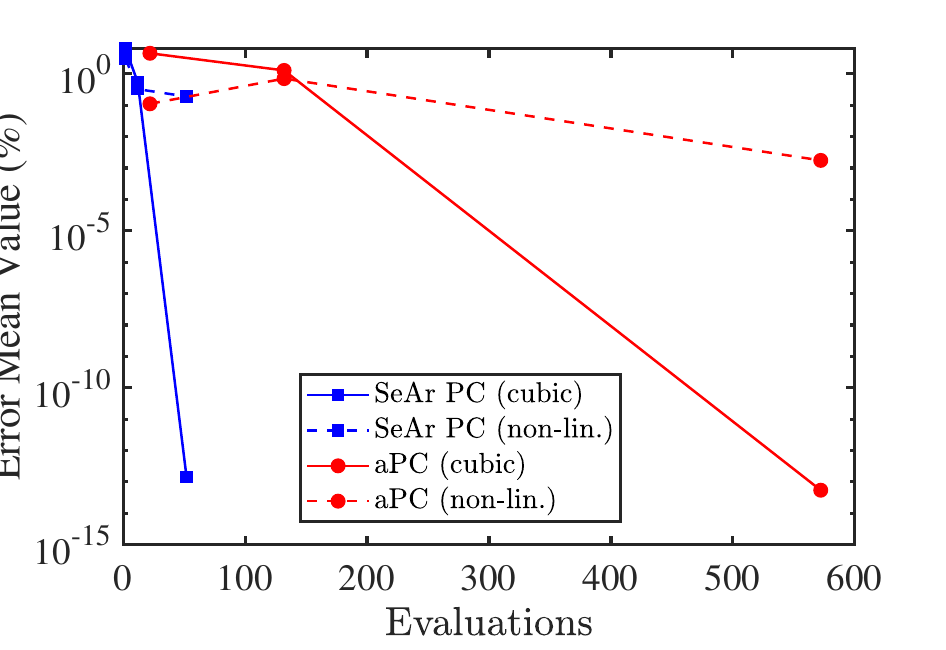}
\includegraphics[width=0.49\textwidth]{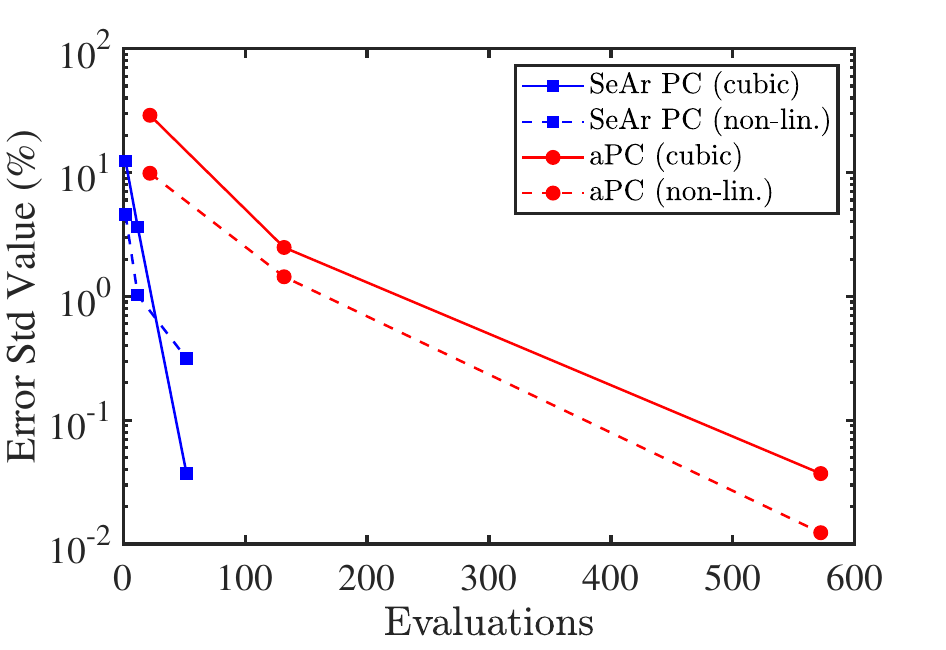}
\includegraphics[width=0.49\textwidth]{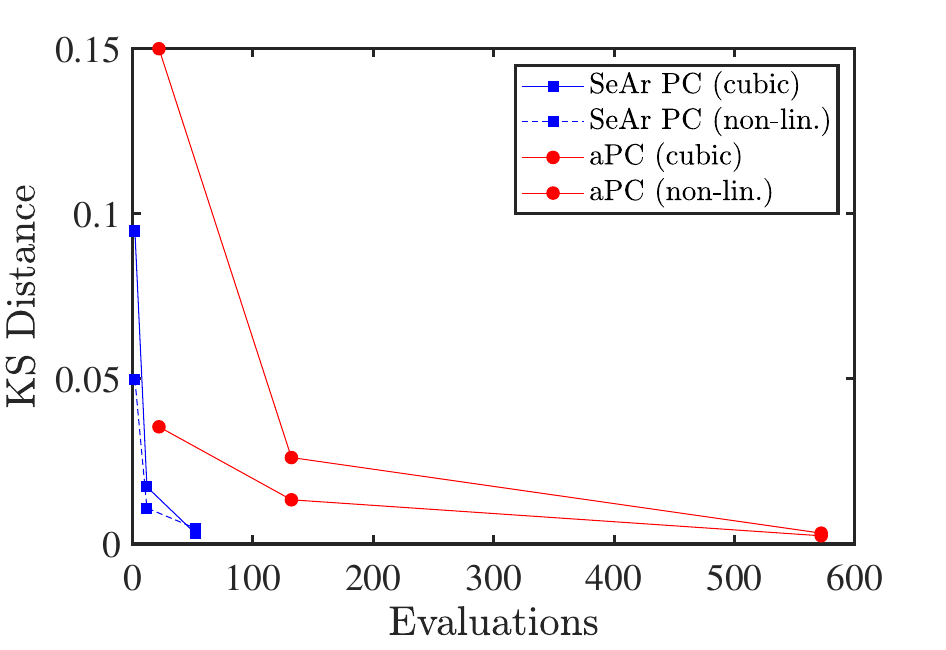}
\end{center}
\caption{Convergence of the first four statistical moments with order/number of evaluations for the multi-modal test case.}
\label{fig:gev} 
\end{figure}

As in the previous case, the histogram of uni-variate samples is plotted in the left panel of Figure \ref{fig:test_inps_gev}, with the histogram for $\mathcal{M}(\boldsymbol{\xi})$ from Monte Carlo sampling on the top right and bottom panels, together with the converging density estimates for SeAr PC. Figure \ref{fig:gev} illustrates the convergence of the first two statistical moments for the various PCE methods. The bottom panel indicates the convergence of the estimated distributions to the direct Monte Carlo samples, as quantified by the KS distance. Again we note that all three methods converged to the true distribution for $f(\mathcal{M}(\boldsymbol{\xi}))$ by $p=3$, however, the computational cost of the PCE with sensitivity enhancement was roughly 10 times less than the standard aPC formulation.  


\begin{table}[!htbp]
\centering
\caption{Convergence of the first four statistical moments with order/number of evaluations for the generalized extreme value test case.}
\begin{tabular}{*1c | *6c}
\toprule
Test case & Method & $p$ & $q$ & $\Delta \mu \,(\%)$ & $\Delta \sigma \,(\%)$ & KS dist. \\
\midrule
\multirow{6}{*}{Cubic function} &
\multirow{3}{*}{SeAr PC} &
1 &2 & 6.34 &12.39 & 0.095 \\
& & 2& 12& 0.54 & 3.62 & 0.017 \\
& & 3& 52& 1.45e-13 & 0.04 & 0.003\\
& \multirow{3}{*}{aPC (WLSQ)}  &
 1 & 22& 4.54 &29.00 & 0.150\\
& &2 &132 & 1.29 & 2.48 & 0.026\\
&&3 &572 & 5.43e-14 & 0.04 & 0.003\\
\midrule
\multirow{6}{*}{Non-lin. function} &
\multirow{3}{*}{SeAr PC} &
1 &2 & 3.08 & 4.57 & 0.050\\
& & 2&12 & 0.33 & 1.02 & 0.011\\
& & 3& 52& 0.19 & 0.32 & 0.005\\
& \multirow{3}{*}{aPC (WLSQ)}  &
 1 &22 & 0.11 & 9.87 & 0.036\\
& &2 &132& 0.71 & 1.44 & 0.013 \\
&&3 &572 & 1.80e-13 & 0.01 & 0.003 \\
\bottomrule
\end{tabular}
\end{table}

\newpage

\subsection{Uncertainty Quantification for a Finite Element Structure determined by Topology Optimisation}
Performing Uncertainty Quantification for a structure designed using Topology Optimisation is challenging due to the size of the uncertain input space involved, with even first-order PCEs requiring more sample evaluations than is practical. In this section we apply SeAr-PC to estimate the effect of parametric uncertainties on the compliance of a two-dimensional Finite Element structure. This structure is an MBB beam, a classic problem in Topology Optimisation (see Sigmund \cite{sigmund200199}). 

\begin{figure}
\begin{center}
\includegraphics[trim={5cm 5cm 5cm 5cm}, width=0.49\textwidth]{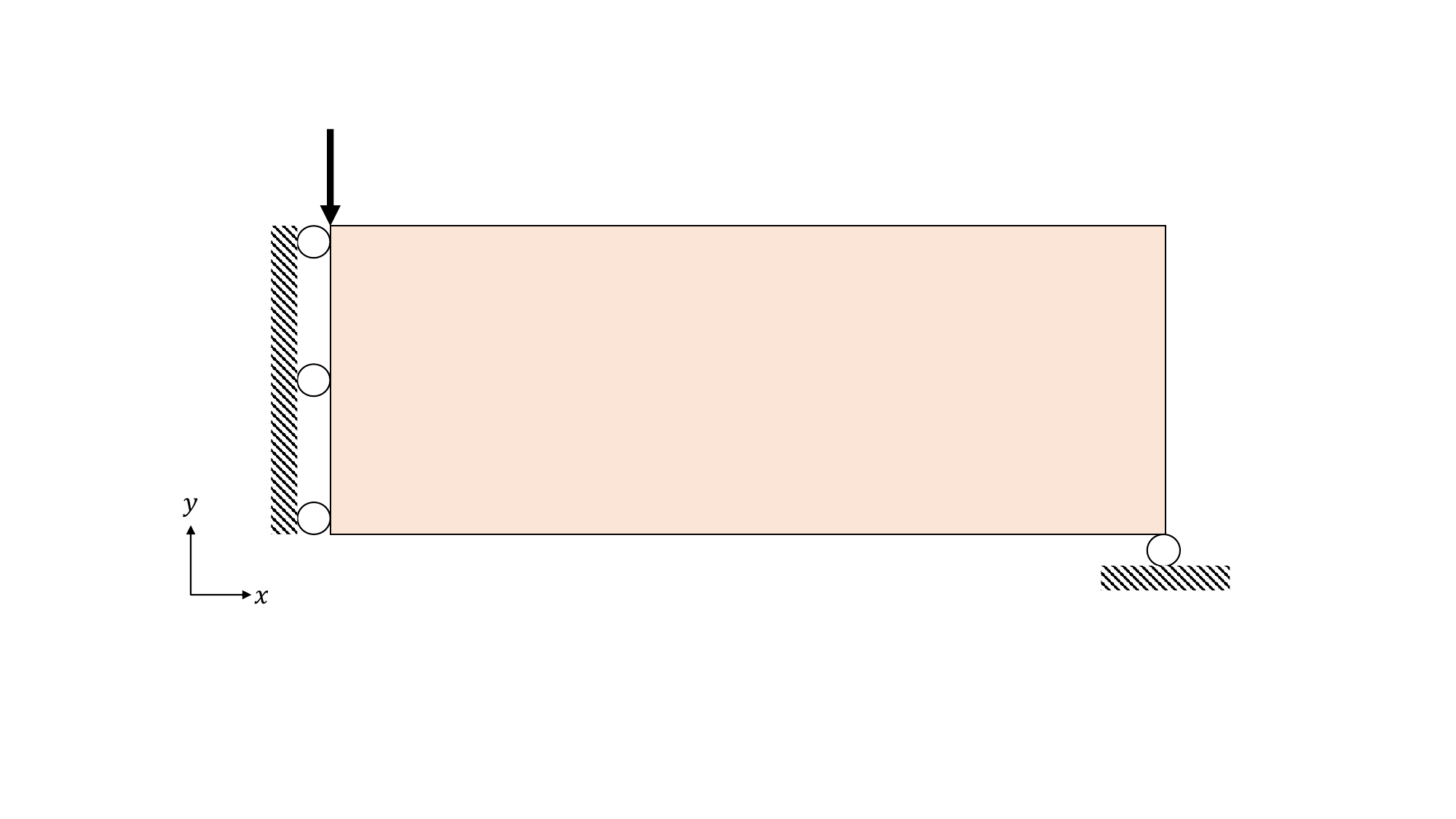}
\includegraphics[trim={0cm 2cm 0 0}, width=0.49\textwidth]{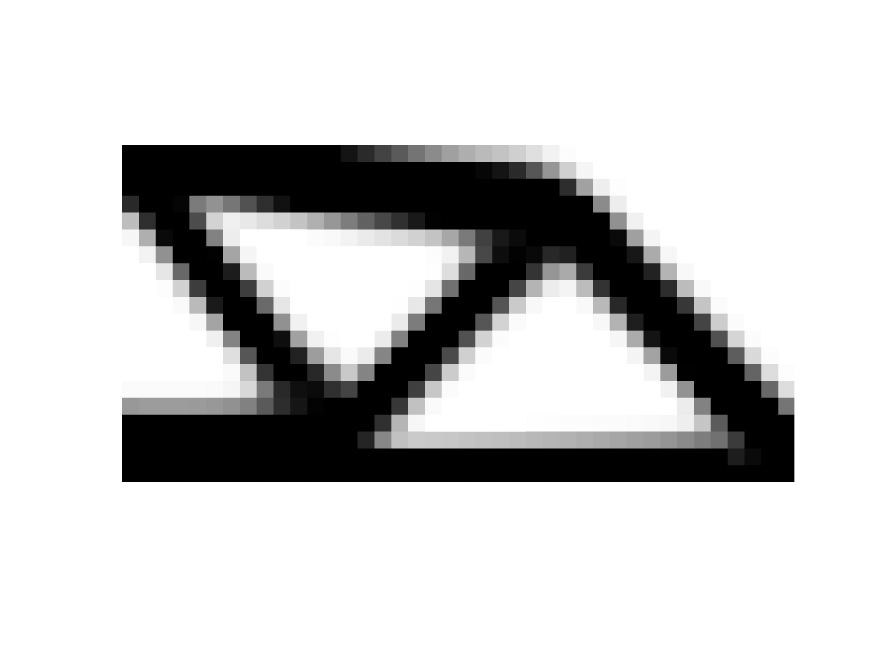}
\end{center}
\caption{Loading and boundary conditions for the MBB beam (left) and optimised material distribution (right).}
\label{fig:struct} 
\end{figure}

Given $n_e=n_x\times n_y$ elements, the SIMP approach is employed to estimate the material distribution, $\boldsymbol{x}\in\Re^{n_x\times n_y}$, where the elements of $\boldsymbol{x}$ determine the densities of the elements within the structure and are constrained to lie in the range $[0,1]$. The loadings and boundary conditions of the MBB beam are illustrated in the left panel of Figure \ref{fig:struct}. The objective of TO is to minimise the compliance of the structure, with a constraint on the minimum amount of material required. Mathematically, the problem is formulated as: 

\begin{align}
    &\underset{\boldsymbol{x}}{\textrm{min}} \;
    \mathcal{C}(\boldsymbol{x})=\boldsymbol{u}^\top K\boldsymbol{u}=\sum_{e=1}^{n_e}E_e(\boldsymbol{x}_e)\boldsymbol{u}_e^\top \boldsymbol{k}_0\boldsymbol{u}_e,
    \label{eq:obj_fun}\\ 
    &\text{such that:}  \begin{cases}
    V(\boldsymbol{x})/V_0=V_f \\
    \boldsymbol{f}=K\boldsymbol{u} \\
    \boldsymbol{0}\leq \boldsymbol{x}\leq \boldsymbol{1}
   \end{cases}    \nonumber
\end{align}

where $\mathcal{C}$ represents the compliance; $\boldsymbol{u}$ and $\boldsymbol{f}$ the global displacement and force vectors; $K$ the global stiffness matrix; $\boldsymbol{u}_e$ the element displacement vector; $k_0$ the element stiffness matrix (assuming $E_e=1)$; $V$ the material volume; $V_0$ the design domain volume; and $V_f$ the prescribed volume fraction. The Young's modulus of the $e$\textsuperscript{th} element, $E_e$, is determined through:

\begin{align}
    E_e(\boldsymbol{x}_e)=E_{\text{min}}+\boldsymbol{x}^3_e(E_0-E_{\text{min}}),
    \label{eq:E}
\end{align}
where $\boldsymbol{x}_e$ is raised to the power three to force it towards either 0 (void) or 1 (solid) and $E_{\text{min}}$ is a small constant to ensure numerical stability. Further details of the optimisation procedure followed to determine the optimal material distribution, $\hat{\boldsymbol{x}}$, and a description of the open source MATLAB code used may be found in Andreassen et al \cite{andreassen2011efficient}. The right panel of Figure \ref{fig:struct} illustrates the optimal material distribution for the problem, using $40\times 20$ elements. 

\begin{figure}
\begin{center}
\includegraphics[width=0.49\textwidth]{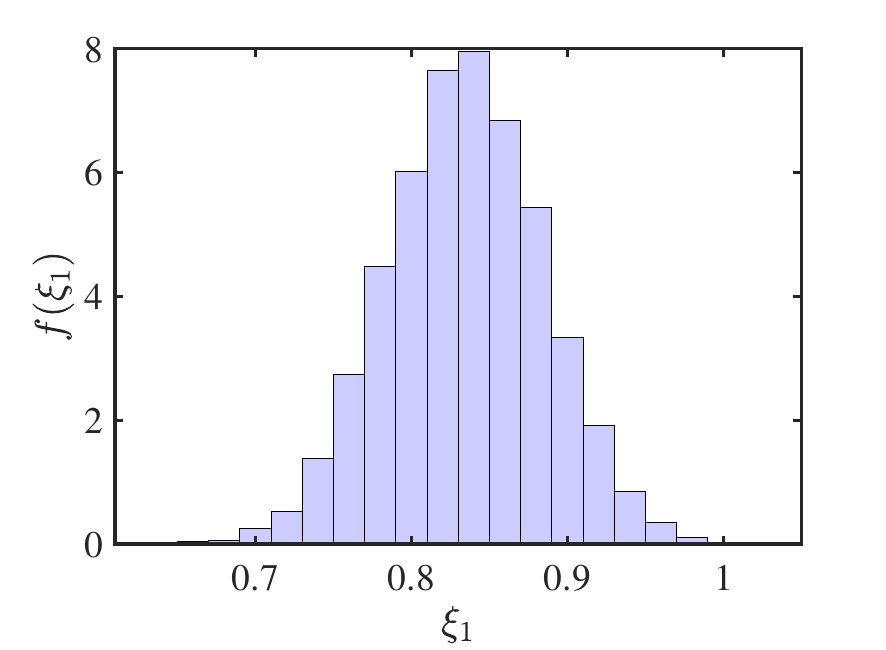}
\includegraphics[width=0.49\textwidth]{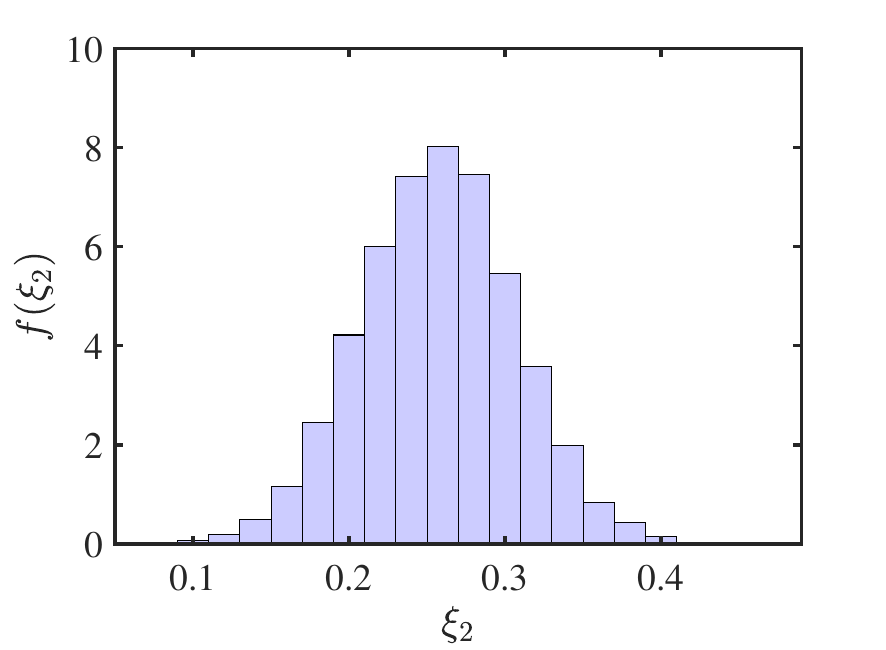}
\includegraphics[width=0.49\textwidth]{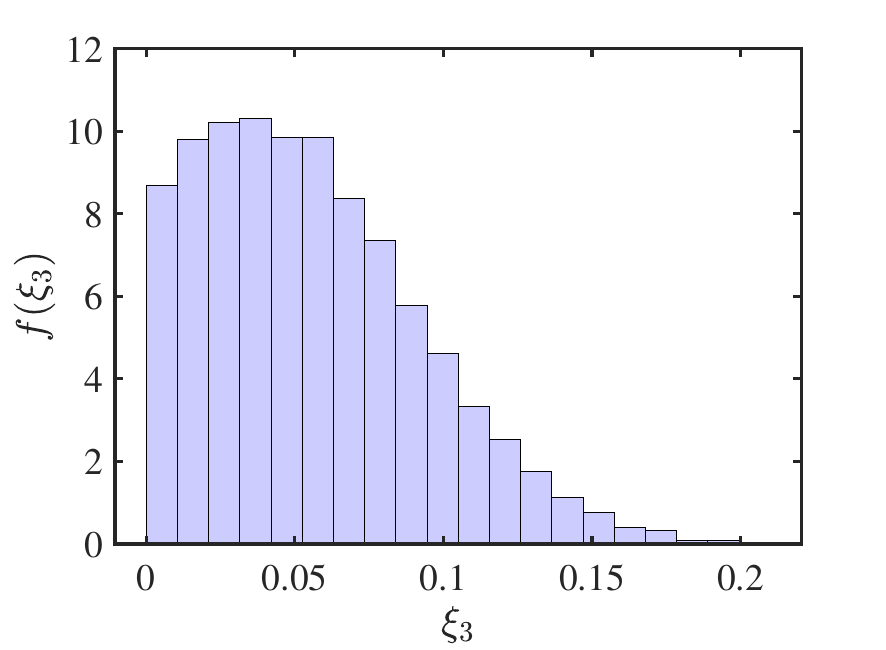}
\includegraphics[width=0.49\textwidth]{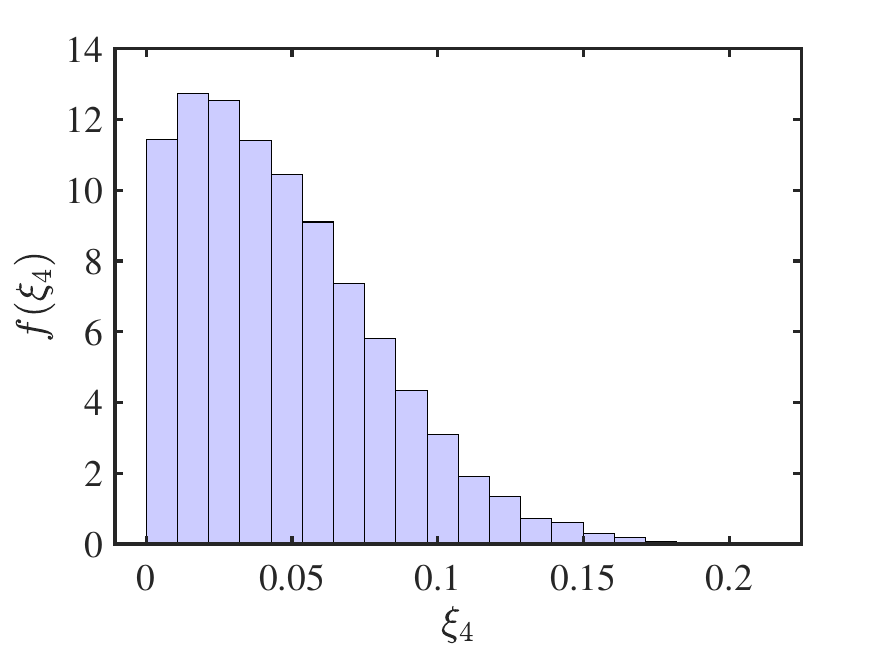}
\end{center}
\caption{Histograms displaying uni-variate Monte Carlo samples from the first four uncertain inputs.}
\label{fig:struct_inps} 
\end{figure}

While the penalty factor in equation \eqref{eq:E} pushes the material distribution towards binary values, 306 elements take intermediate values. The density for these uncertain elements is denoted $\hat{\boldsymbol{x}}_s$. In practice manufacturing imperfections will cause over- or under-deposition for these elements. In this section we investigate the effect of this uncertainty on the compliance of the structure, using the assumed joint density of $N(\hat{\boldsymbol{x}}_s, 0.05\boldsymbol{1}, \Sigma)$ for these 306 elements, where $\Sigma$ is equal to the identifity matrix. The joint density is truncated at $\boldsymbol{0}$ and $\boldsymbol{1}$ to reflect the bounds on $\boldsymbol{x}$. The joint density for the first four uncertain parameters is plotted in Figure \ref{fig:struct_inps}. As can be seen from the bottom two panels, the truncation can have a significant effect on the normality of the probability distribution. 

\begin{figure}
\begin{center}
\includegraphics[width=0.49\textwidth]{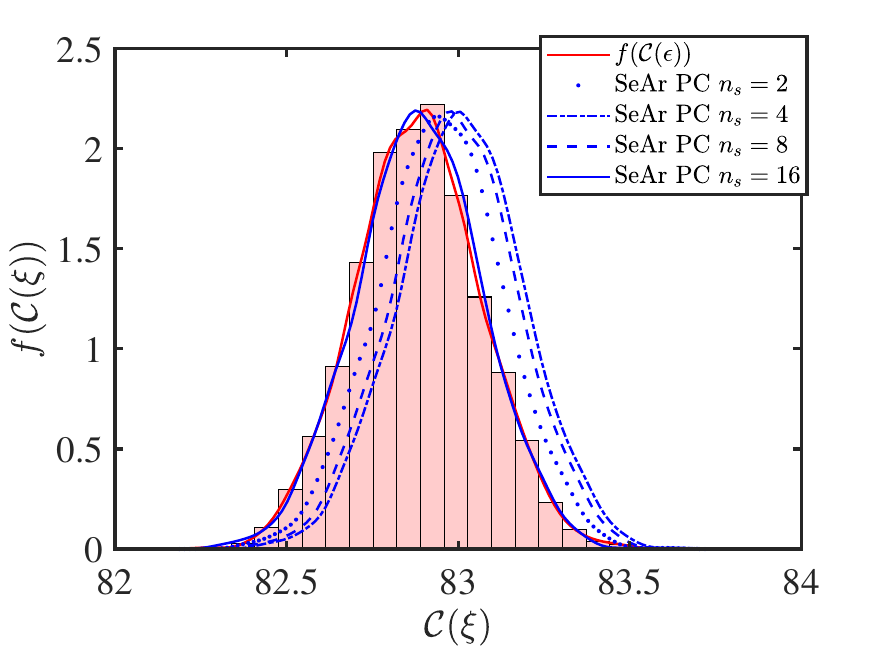}
\end{center}
\caption{Comparison between a probability distribution for the compliance using 10,000 Monte Carlo samples (red) and first-order SeAr-PC (blue) for varying values of the oversampling ratio.}
\label{fig:struct_out} 
\end{figure}

The response of the compliance to uncertain deposition for the elements in $\hat{\boldsymbol{x}}_s$ is represented using a PCE. A single realisation of $\mathcal{C}$ requires the solution of the matrix equation $F=KU$, which will be relatively expensive for a large structure. With such a high dimensional uncertain space, more samples will be required to estimate the PCE coefficients with the WLSQ and Smolyak approaches than is practical to run, even if only first-order PCEs are used. However, recall that for a first order PCE with sensitivity enhancement, the required number of samples is decoupled from the dimensions of the uncertain space. With sensitivity enhancement, a number of forwards evaluations equal to the oversample ratio is required. Such a formulation is particularly convenient for this problem as the compliance is easily differentiated with respect to $\boldsymbol{x}$, leading to analytical expressions for the sensitivity:

\begin{align}
    \frac{\partial \mathcal{C}}{\partial \hat{\boldsymbol{x}}_e}=-3\hat{\boldsymbol{x}}_e^2(E_0-E_{\text{min}})\boldsymbol{u}_e^\top \boldsymbol{k}_0\boldsymbol{u}_e, 
\end{align}
the block equations for the sensitivities can therefore be obtained for negligible cost. Figure \ref{fig:struct_out} compares the histograms for $\mathcal{C}$ estimated using a first-order PCE with one obtained from 10,000 Monte Carlo samples. In Table \ref{tab:TO} we tabulate the first two statistical moments of the probability distribution for $\mathcal{C}$, comparing SeAr-PC against aPC with WLSQ for varying values of the oversampling ratio. The table demonstrates the benefit of the sensitivity enhancement when scaling to very high-dimensional uncertain spaces: convergence can be achieved for a handful of function evaluations using the sensitivity enhancement, instead of the hundreds required by WLSQ. Figure \ref{fig:struct_out} illustrates the matching between the Monte Carlo estimate for $\mathcal{C}$ between the Monte Carlo sampling (red) and a Gaussian distribution with the statistical moments determined by SeAr-PC with an oversampling ratio of 8. To minimise sampling error, the statistical moments are determined analytically from the PCe coefficients, as in \eqref{eq:pce_mom}.

\begin{table}
\begin{center}
\caption{First two statistical moments of $\mathcal{C}$, estimated by Monte Carlo sampling, SeAr-PC and aPC with WLSQ, with varying values for the oversampling ratio.}
\begin{tabular}{*4c}
\toprule
Method & $q$ & $\mu$  & $\sigma$  \\
\midrule
Monte Carlo sampling & 10,000  & 82.885 & 0.182 \\
\midrule
\multirow{5}{*}{SeAr-PC ($p=1$)} 
&1 & 82.805 & 0.182 \\
&2 & 82.948 & 0.184 \\
&4 & 83.009 & 0.184 \\
&8 & 82.987 & 0.184 \\
&16 & 82.887 & 0.183 \\
\midrule
\multirow{2}{*}{aPC ($p=1$)}
& 307 & 82.703 & 2.461 \\
& 614 & 82.879 & 0.1720\\

\bottomrule
\label{tab:TO}

\end{tabular}
\end{center}
\end{table}

\begin{figure}
\begin{center}
\includegraphics[trim={0cm 2cm 0 0},width=0.8\textwidth]{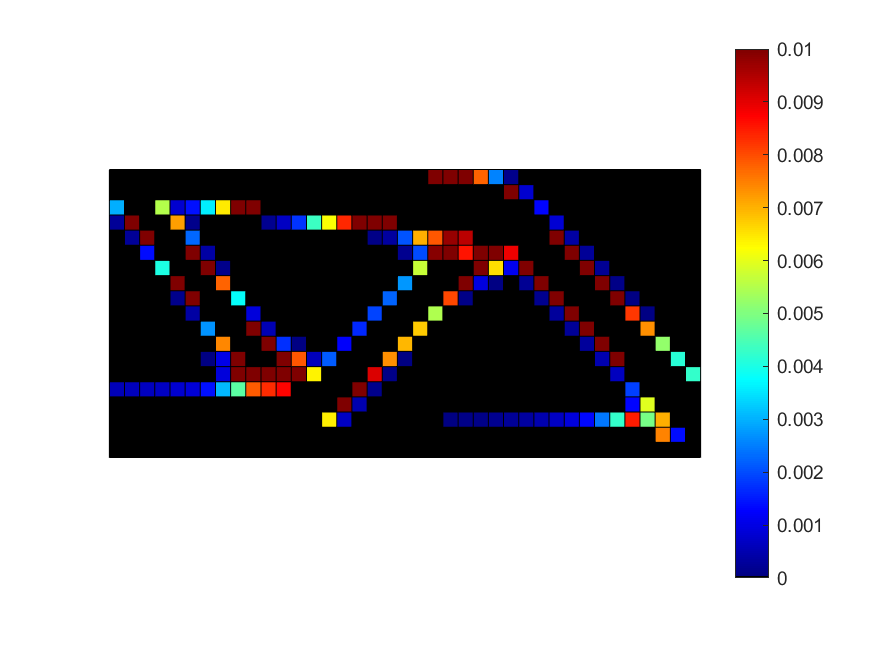}
\includegraphics[trim={0cm 2cm 0 0},width=0.8\textwidth]{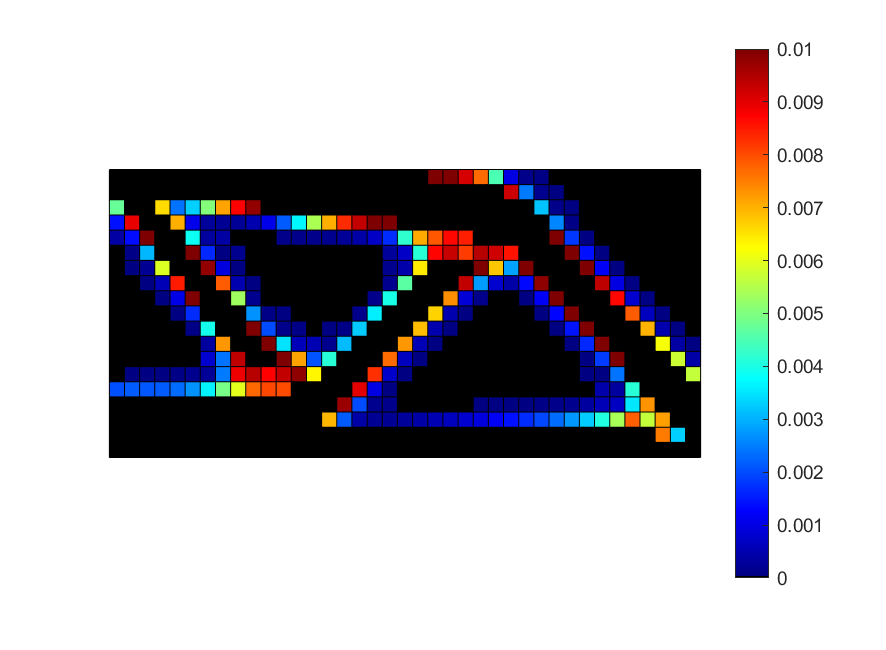}

\end{center}
\caption{Contributions of each element to the total uncertainty, quantified by the first-order Sobol indices (top) and mean sensitivity of the MC samples (bottom). }
\label{fig:struct_sobol} 
\end{figure}

A desirable feature of PCEs is that, having estimated the coefficients $\boldsymbol{\lambda}$, a first-order Sobol sensitivity analysis can be performed at negligible cost, quantifying the relative contribution of each uncertain parameter to the total uncertainty. The Sobol sensitivity index for the $j$\textsuperscript{th} uncertain input, $S_j$, is estimated through \cite{crestaux2009polynomial,pepper_mmodal}:

\begin{align}
    S_j=\frac{\sum_{j\in I_j}\boldsymbol{\lambda}_j^2 \langle \Psi_j, \Psi_j \rangle}{\sum_{j=1}^P\boldsymbol{\lambda}_j^2 \langle \Psi_j, \Psi_j\rangle},
\end{align}
where the numerator represents those terms in the multi-index in which the $j$\textsuperscript{th} uncertain input is non-zero. The polynomial norms were calculated using Smolyak's algorithm with a level 1 sampling grid. In Figure \ref{fig:struct_sobol} the Sobol indices are presented as a heatmap across the optimised material distribution, allowing the designer to visualise the spatial distribution of the main sources of uncertainty in the design. For validation, this heatmap is compared to one obtained through an analysis of the Monte Carlo samples, with the average sensitivity for the $k$\textsuperscript{th} element computed as:

\begin{align}
    S_k\approx \frac{\sum_{m=1}^{q}\big|\frac{\partial \mathcal{C}}{\partial \boldsymbol{\xi}^{(m)}_k}\big|}{\sum_{i=1}^{n_u}\sum_{m=1}^{q}\big|\frac{\partial \mathcal{C}}{\partial \boldsymbol{\xi}^{(m)}_i}\big|}
\end{align}
As can be seen from the Figure, the two heatmaps strongly resemble one another, with both indicating that the uncertainty is concentrated in the corners of the structure, and on the right-hand edge.

\newpage

\section{Conclusion}
{In this paper the arbitrary Polynomial Chaos formulation is extended through combination with a recently method for enhancing PCEs with sensitivity information. The hybrid method proposed here} is intended to address two common criticisms of PCEs: that the range of stochastic processes for which an optimal polynomial basis exists is limited; and that the required number of function evaluations increases rapidly with the size of the uncertain input space. 

Incorporating sensitivity information addresses the latter concern by obtaining an additional $n_u$ equations for each model evaluation. This sensitivity information could be analytic (as is the case here) or derived through the adjoint formulation for minimal additional cost. For PCEs of order $p\geq2$ the proposed method scales more efficiently with increasing $n_u$: requiring a factor of $n_u$ fewer model evaluations than Smolyak quadrature and the standard LSQ methods. Interestingly, for a first-order PCE the number of required model evaluations is decoupled from $n_u$. This feature opens up the application of PCEs to a wider set of potential applications, which were previously too high-dimensional for even first-order PCEs, as they required many more model applications than were feasible. Furthermore, the use of arbitrary Polynomial Chaos allows optimal basis polynomials to be found for a wider variety of stochastic processes, which do not have an optimal basis in the Askey scheme. 

We demonstrate the two features of the method through applications to synthetic test cases involving both high-dimensional uncertain spaces and a range of stochastic processes: including multi-modal and truncated probability distributions. In particular, we focus on an application to a toy problem in Topology Optimisation, featuring 306 uncertain parameters. For this case, we demonstrate that the proposed method can produce a reasonable estimate of the first two statistical moments of the compliance, requiring only a handful of model evaluations. The test case also demonstrates two advantages of PCEs for this type of application: firstly that sampling error can be minimised by computing the statistical moments directly from the PCE coefficients; and secondly, that a Sobol sensitivity analysis of these coefficients can indicate to an engineer the regions in the structure that have the most contribution to the overall uncertainty. The developments described here unlock a range of potential applications of PCEs to problems in Topology Optimisation, additive manufacturing, and robust design that could be explored further in future work.

\newpage
\bibliography{main}

\end{document}